\documentclass{article}
\usepackage{graphicx} 
\usepackage{caption}
\usepackage{subcaption}

\usepackage[utf8]{inputenc}
\usepackage{fullpage}
\usepackage{amsmath,amssymb,amscd,amsthm,amsfonts}
\usepackage{amssymb}
\usepackage{hyperref}
\usepackage{amsmath}
\usepackage[section]{placeins}
\usepackage{parskip}
\usepackage[colorinlistoftodos,bordercolor=orange,backgroundcolor=orange!20,linecolor=orange,textsize=scriptsize]{todonotes}
\usepackage{amsthm}
\usepackage{verbatim}
\usepackage{thmtools}
\theoremstyle{plain}
\usepackage{pgf,tikz}
\usepackage{dirtytalk}
\usepackage{pdfpages}
\usepackage{mathrsfs}
\usepackage{enumitem}   

\usetikzlibrary{arrows}
\usepackage{authblk}
\usepackage[labelformat=simple]{subcaption}

\usepackage{
hyperref,
}

\definecolor{qqzzqq}{rgb}{0.,0.6,0.}
\definecolor{qqwuqq}{rgb}{0.,0.39215686274509803,0.}
\definecolor{ffqqqq}{rgb}{1.,0.,0.}
\definecolor{ududff}{rgb}{0.30196078431372547,0.30196078431372547,1.}

\newenvironment{Proof1}{\emph {Proof of Lemma 1.}}{\hfill$\square$}

\newcommand\restr[2]{{
  \left.\kern-\nulldelimiterspace 
  #1 
  \littletaller 
  \right|_{#2} 
  }}

\newcommand{\littletaller}{\mathchoice{\vphantom{\big|}}{}{}{}}

\usepackage{subfiles}

\DeclareMathOperator{\EEmbed}{\mathrm{Embed}}

\DeclareMathOperator{\cir}{cr}
\DeclareMathOperator{\ess}{ess}
\DeclareMathOperator{\verti}{vert}

\DeclareMathOperator{\DDiag}{\mathrm{Diag}}
\DeclareMathOperator{\DDiam}{\mathrm{Diam}}
\DeclareMathOperator{\ddiam}{\mathrm{Diam}}

\newtheorem{theorem}{Theorem}
\newtheorem{definition}{Definition}
\newtheorem{conjecture}{Conjecture}

\newtheorem{lemma}{Lemma}
\newtheorem{proposition}{Proposition}
\newtheorem{corollary}{Corollary}

\title{Borsuk and Vázsonyi problems through Reuleaux polyhedra}
\author[1, 2] {Gyivan López-Campos\footnote{Partially supported by grant PAPIIT IG100721}}
\author[1]{Déborah Oliveros\footnote{Partially supported by grant PAPIIT IG100721}}
\author[2]{Jorge L. Ramírez Alfonsín\footnote{Partially supported byAPP INSMI-CNRS}}
\affil[1]{Universidad Nacional Aut\'onoma de México, Instituto de Matemáticas, Unidad Juriquilla}
\affil[2]{IMAG, Univ. Montpellier, CNRS, Montpellier, France}
\date{\today}

\begin{document}

\maketitle
\begin{abstract}
    The Borsuk conjecture and the V\'azsonyi problem are two attractive and famous questions in discrete and combinatorial geometry, both based on the notion of diameter of bounded sets. In this paper, we present an equivalence between the critical sets with Borsuk number 4 in $\mathbb{R}^3$ and the minimal structures for the V\'azsonyi problem by using the well-known Reuleaux polyhedra. The latter leads to a full characterization of all finite sets in $\mathbb{R}^3$ with Borsuk number 4.

    The proof of such equivalence needs various ingredients, in particular, we proved a conjecture dealing with \emph{strongly critical configuration} for the Vázsonyi problem and showed that the diameter graph arising from involutive polyhedra is vertex (and edge) 4-critical.
\end{abstract}
\section{Introduction}
 \emph{The Borsuk partition} and \emph{The frequent large distance problems} are two attractive and well-known questions in discrete and combinatorial geometry, both based on the notion of \emph{diameter} of bounded sets. The \emph{diameter} of a bounded set $S\subset \mathbb{R}^d$ is defined as $\ddiam (S):= \sup\limits_{x,y \in S}\{\|x-y\|\}$ where $\|x,y\|$ denotes the Euclidean distance between $x$ and $y$. If $S$ is a finite set of points, the diameter is the maximum Euclidean distance between any two points of $S$. In this paper we put forward an equivalence of these problems by considering their finite \emph{strongly critical} configurations.

 In 1933, Borsuk \cite{borsuk1933drei} proposed the following question (sometimes known as \emph{Borsuk conjecture}) :
  \begin{quote}
      Is every set $S\subset \mathbb{R}^d$ with finite diameter $\DDiam (S)$ the union of at must $d+1$ sets of diameter less than $\DDiam (S)$?
  \end{quote}

It is known to be true for $d=2$ (see \cite{borsuk1933drei}) and for $d=3$ (see \cite{perkal1947subdivision}, \cite{eggleston1955covering} and \cite{grunbaum1957simple} for a simpler proof).
  
For fifty years, Borsuk's conjecture was believed to be true until 1993 when Kahn and Kalai \cite{kahn1993counterexample} proved it to be false for $d=1325$ and for each $d>2014$. Nowadays, there are known counterexamples in dimensions 64 and higher \cite{jenrich201464} but the problem is still open for $4\leq d\leq 63$. We refer the reader to \cite{rauigorodskiui2004borsuk} for a survey on the Borsuk conjecture.

Recall that the \emph{Borsuk number} of a bounded set $S\subset \mathbb {R}^d$, denoted by $a(S)$, is the smallest number of subsets that $S$ can be partitioned, such that each part has smaller diameter than $S$. Also, recall that the \emph{diameter graph} $\DDiam_V$ of finite $V\subset \mathbb{R}^3$ is defined as the graph with set of vertices $V$ and two vertices are joined by an edge if their distance is a diameter. These are helpful definitions in order to deal with the Borsuk problem for a finite set of points $V$, since in this case the equality $\chi(\DDiam_V)=a(V)$ holds, where $\chi (G)$ denotes the chromatic number of the graph $G$. 

 Boltyanski proved that in $\mathbb{R}^2$ a bounded set is not the union of two sets with smaller diameter if and only if it has a unique \emph{completion} to a body of constant width see \cite{Bol} for the original proof in Russian (or \cite[pp-245]{boltyanski2012excursions} in English).

By using the above definitions, Boltyanski characterized all the sets in $\mathbb{R}^2$ having Borsuk number 3 (that is, attaining the maximum). Unfortunately, the same argument does not work in $\mathbb{R}^3$ for the sets with Borsuk number 4. For instance, the set of vertices of a regular tetrahedron has Borsuk number 4 but its completion to a body of constant width is not unique (see \cite{meissner1918durch}, \cite{montejano2017meissner}, \cite[pp-358]{martini2019bodies}). 

In the same spirit, it turned out to be a challenging problem to characterize all the sets in $\mathbb{R}^3$ with Borsuk number 4. In \cite{hujter2014multiple}, Hujter and Lángi give all the configurations of these sets up to 7 points and aforementioned,  we cite:
\begin{quote}
    \say{A complete characterization of the Borsuk number of finite sets in $\mathbb{R}^3$, even of those with $a(S)=4$, looks hopeless.}
\end{quote}

Our main result gives a complete characterization of finite subsets in $\mathbb{R}^3$ with Borsuk number 4. We do so by using some recent tools/results about involutive polyhedra and by characterizing the \emph{critical} Borsuk configurations, that is, the finite sets not having subsets with the same Borsuk number.

Our approach is closely related to the \emph{frequent large distance problem}:

\begin{quote}
     Given $0<d<n$, what is the maximum number of diameters over all the sets of $n$ points in $\mathbb {R}^d$?
\end{quote}

We denote by $e(d,n)$ such maximum number of diameters. This is one of the oldest problems in discrete and combinatorial geometry. It was first proposed in 1934 by Hopf and Pannwitz \cite{hopf1934aufgabe} in the plane and then generalized to all dimensions. 

Given a finite set $V\subset \mathbb {R}^d$, we let $e(V)$ be the number of diameters in $V$ (we keep the same notation introduced in \cite{kupitz2010ball}). We say that $V$ is an \emph{extremal configuration} for the frequent large distance problem if $e(V)=e(d,|V|)$.

It is well known that $e(2,n)=n$ and how all the extremal configurations look like (see \cite[pp 213-214]{pach2011combinatorial},  \cite{kupitz1979extremal}). For $d=3$, the problem is better known as \emph{the Vázsonyi problem} in honor to Vázsonyi, who conjectured that $e(3,n)=2n-2$. Grünbaum \cite{grunbaum1956proof}, Heppes \cite{heppes1956beweis} and Straszewicz \cite{straszewicz1957probleme} proved independently to be true and more recently, a completely different proof was given by K. Swanepoel \cite{swanepoel2008new}.  Finally,  Kupitz, Martini and Perles \cite{kupitz2010ball} characterize all the extremal configurations.

We say that $V$ is a \emph{critical configuration} for the Vázsonyi problem if $V$ is an extremal configuration and any point 
of $V$ is adjacent to at least 3 diameters. We also say that $V$ is \emph{strongly critical} if $V$ does not have an extremal configuration subset. By using the characterization of the extremal configurations, we have that being strongly critical implies to be critical, however the opposite direction is not true.

The existence of a set of 8 points that is critical but not strongly critical was claimed in \cite{kupitz2010ball} and intended to be given in a future work, however,  as far as we are aware it was never published. By using bodies of constant width, we are able to construct an explicit critical configuration of 8 points in $\mathbb {R}^3$ that is not strongly critical (see end of Section \ref{ReuleauxSection}).

Our approach led us to investigate the \emph{ball polyhedra} (that is, ball polytopes in dimension 3).  In \cite{bezdek2007ball}, it was proved that the 1-skeleton of ball polyhedra arising from extremal set of points in $\mathbb {R}^3$ is a 2-connected planar graph (this was already observed in \cite{bezdek2007ball}).  In the same paper, the authors also proposed the following\\

\begin{conjecture}\cite{kupitz2010ball}\label{vazcritical}
    Let $V\subset\mathbb{R}^3$ be an extremal set.  Then, $\mathcal{B}(V)$ is a standard ball polyhedron if and only if $V$ is strongly critical. 
\end{conjecture}

We were able to prove this conjecture in (Lemma \ref{3connected}). The latter yields to a nice equivalence between strongly critical configurations for the Vázsonyi problem and the \emph{Reuleaux polyhedra} (Theorem \ref{Conjballteo}). Furthermore,  the 4-criticality of the \emph{diagonal} graph arising from \emph{involutive polyhedra} (Lemma \ref{critical}), led us to a full characterization of all finite sets in $\mathbb{R}^3$ with Borsuk number 4.\\

\begin{theorem}\label{MainTheorem}
Let $V\subset \mathbb{R}^3$ be a finite set with unit diameter and $\mid V \mid =n \geq 4$. The following statements are equivalent
    \renewcommand{\theenumi}{\roman{enumi}}
    \begin{enumerate}
        \item $V$ has a subset that is an extremal configuration for the Vázsonyi problem.
        \item $V$ has Borsuk number 4.
        \item There is a $V_1\subset V$ such that $\mathcal{B}(V_1)$ is a Reuleaux polyhedron.
    \end{enumerate}
\end{theorem}

The organization of the paper is as follows. In the next section we present a number of results and notions needed for the rest of the paper. In particular, we discuss some background on both the \emph{ball polyhedra} and   \emph{Reuleaux polyhedra} as well as their properties. In Section \ref{KeyLemma}, we prove a key lemma on the chromatic number of the \emph{diagonal} graph of involutive polyhedra. This is not only interesting for its own sake, but it is a crucial brick for our contributions. Section \ref{Main Results} is mainly  devoted to prove our main results. We finally end with some concluding remarks.

\section{Preliminaries}

We review some results and notions on ball polyhedra and Reuleaux polyhedra needed throughout the paper.  We refer the reader to \cite[pp 132-141]{martini2019bodies} for further details. We also discuss some useful background on involutive polyhedra. 

\subsection{Ball-polyhedra}\label{ball}

Given a finite subset $V$ of $\mathbb{R}^3$, the ball set of $V$ is defined as $\mathcal{B}(V)=\{y\in \mathbb{R}^3: \forall x \in V, \|x-y \|\leq 1 \}$. If the radius of the \emph{circumball} of $V$, denoted by $\cir(V)$, is less than 1, then $\mathcal{B}(V)$ is called the \emph{ball polyhedron} associated with $V$. A point $v\in V$ is \emph{essential} if $\mathcal{B}(V)\subsetneqq  \mathcal{B}(V\backslash \{v\})$. The subset of essential points will be denoted as $\ess (V)$. A finite set $V\subset \mathbb{R}^3$ satisfying $\cir(V)<1$ and $V=\ess(V)$ is  called \emph{tight}.

The following four theorems are due to Martini, Kupitz and Perles \cite{kupitz2010ball}.\\

\begin{theorem}\label{eqtight}\cite{kupitz2010ball}
    Assume that $V\subset \mathbb{R}^3$ is finite and $\ddiam V=1$. Then 
    \begin{enumerate}
        \item $\cir (V)<1$
        \item If a point $v\in V$ is incident with (at least) two diameters of $V$, then $v\in \ess(V)$.
        \item If $V$ is extremal for the Vázsonyi problem, then $V$ is tight.
    \end{enumerate}
\end{theorem}
\begin{definition}\label{vertex}
Let $V$ be a tight set of points.  The \emph{facial structure} of the ball polyhedron $\mathcal{B}(V)$ is defined as
    \begin{enumerate}

        \item For a essential point $p\in V$ the set $F_p:=\{x\in \mathcal{B}(V): ||x-p||=1\}$ is a \emph{facet} of $\mathcal{B}(V)$.
        
        \item     A boundary point $z$ of $\mathcal{B} (V)$ is a \emph{vertex} of $\mathcal{B}(V)$ if either $z$ belongs to three or more distinct facets of $\mathcal{B} (V)$, in this case $z$ is called a \emph{principal vertex}, or $z\in V\cap \mathcal{B} (V)$ and $z$ belongs to exactly two facets of $\mathcal{B}(V)$, in this case $z$ is called a \emph{dangling vertex}. Denote by $\verti \mathcal{B}(V)$ the set of vertices of $\mathcal{B}(V)$. In other words, $z\in \verti \mathcal{B}(V)$ if and only if $z \in \mathcal{B}(V)$ and $\|z-p\|=1$ holds for at least three points $p\in V$, or if $z\in V\cap \mathcal{B} (V)$ and $\|z-p\|=1$ holds for exactly two points $p\in V$.
        
        \item An \emph{edge} of $\mathcal{B}(V)$ is the closure of a connected component of $(F_p\cap F_q)\backslash (\verti \mathcal{B} (V))$, where $\{p, q\}$ ranges over all pairs of distinct points of $V$.
        
        \item The set of faces of $\mathcal{B}(V)$, including facets, edges, vertices and improper faces $\mathcal{B}(V)$ and $\emptyset$, is the \emph{spherical face complex} of $\mathcal{B} (V)$ denoted by $\mathcal{SF(B}(V))$. In particular, the 1-skeleton of $\mathcal{SF(B}(V))$ is the set of vertices and edges of $\mathcal{B}(V)$ viewed as a graph. 
    \end{enumerate}
\end{definition}

\begin{theorem}\label{2connected}\cite{bezdek2007ball}
Given a tight finite set $V\subset \mathbb{R}^3$ and $|V|\geq 3$, the 1-skeleton of $\mathcal{SF(B}(V))$ is planar and 2-connected. 
\end{theorem}

The converse of the above theorem is shown in \cite{almohammad2020analogue}.
The following result was called the \emph{extended GHS Theorem} in \cite{kupitz2010ball} after Gr\"umbaum, Heppes and Straszewicz who gave the proofs for the Vázsonyi problem independently.

\begin{theorem}\label{GHS}\cite{kupitz2010ball}
    \textbf{(GHS)} Let $V\subset \mathbb{R}^3$ be finite with $\mid V \mid =n \geq 4$ and $\ddiam V=1$. The following three statements are equivalent
    \begin{enumerate}
        \item V is extremal for the Vázsonyi problem, i.e., $e(V)=e(3, n)$.
        \item $e(V)=2n-2$.
        \item $V$ is tight and $V=\verti \mathcal{B} (V)$.
    \end{enumerate}
\end{theorem}

An \emph{involutory self-duality} of $\mathcal{SF(B}(V))$ is an order reversing map $\varphi: \mathcal{SF(B}(V))\rightarrow{\mathcal{SF(B}(V))}$ of order two ($\varphi^2=Id$) and that sends every vertex $v\in \mathcal{SF(B}(V))$ to its corresponding \emph{ dual face} $F_v\in \mathcal{SF(B}(V))$. This involution can be naturally extended to the edges as follows: for every edge $ab\in \mathcal{SF(B}(V))$, $\varphi(ab)=\varphi(a)\varphi(b)$ is the edge induced by the intersection of $F_a$ and $F_b$.

\begin{theorem}\label{dualedges} \cite{kupitz2010ball}
Let $V$ be an extremal Vázsonyi configuration in $\mathbb{R}^3$. Then, there is always a unique edge-induced involution $\varphi: \mathcal{SF(B}(V))\rightarrow{\mathcal{SF(B}(V))}$ such that for all $v\in V$,  $v\notin \varphi (v)$.
\end{theorem}

We will refer to this involution as \emph{the canonical involution}.

A ball polyhedra $\mathcal{B}(V)$ is called \emph{standard} if $\mathcal{SF(B}(V))$ is a \emph{polytopal} lattice (that is, $\mathcal{SF(B}(V))$ can be realized as the face lattice of a 3-polytope). Numerous papers have focused their attention in studying this kind of ball polyhedra. For instance, it is known that $Q$ is a standard ball polyhedron if and only if either for any supporting sphere $\mathbb {S^d}(p,r)$ of $Q$, the intersection $Q\cap\mathbb {S^d}(p,r) $ is homeomorphic to a closed Euclidean ball of some dimension \cite{bezdek2007ball} or the intersection of two faces is either empty, a vertex or an edge \cite{montejano2020graphs} (see also  \cite{bezdek2007ball, montejano2017meissner}).

In \cite{kupitz2010ball}, it was aforementioned that not all the extremal configurations for the Vázsonyi problem induce a standard ball polyhedron. The example that we present in Section \ref{SpecialConfig} is a critical configuration for the Vázsonyi problem, but it turns out not to be a standard ball polyhedron  (see Figure \ref{Example8}).


\subsection{Reuleaux polyhedra}\label{Graphs}

A standard ball polyhedron $\mathcal{B}(V)$ satisfying $V=\verti\mathcal{B}(V)$ is called a \emph{Reuleaux} polyhedron, and denoted by $R(V)$. Reuleaux polyhedra enjoy several attractive properties. For instance, they are ``frames" of bodies of \emph{constant width} in $\mathbb{R}^3$; see for example, the \emph{Meissner} polyhedra constructed in \cite{montejano2017meissner} or the \emph{Peabodies} built in \cite{arelio2022peabodies}.

It is known that the set of vertices of a Reuleaux polyhedron $V$ form an extremal configuration for the Vázsonyi problem \cite{montejano2020graphs}. Furthermore, by using the density of the Reuleaux polyhedra in the set of bodies of constant width (investigated in \cite{sallee1970reuleaux}), it was 
showed in \cite{hujter2014multiple}, that the vertex set of a Reuleaux polytope has Borsuk number 4.  This fact can also be deduced from \cite[Theorem 3]{montejano2020graphs} where the chromatic number for the diameter graph of $V$ was shown to be equal 4.

A graph $G$ is called \emph{polyhedral graph} if it is a simple, 3-connected, planar graph. The name comes after Steinitz' characterization \cite{steinitz1922polyeder} stating that $G$ is a polyhedral graph if and only if it is the 1-skeleton of a convex 3-polytope. Since the Reuleaux polyhedra are standard ball polyhedra, then they have polytopal structure and hence their 1-skeleton is a polyhedral graph.

\subsection{Involutive graphs}\label{Invgraphs}

Let $G=(V,E)$ be a planar graph with set of vertices $V$ and set of edges $E$ and let $G^*$ be its \emph{dual} graph (we refer the reader to \cite{bondy2008graph} for graph theory notions).  We say that $G$ is {\em self-dual} if $G$ is isomorphic to $G^*$, that is, there is a bijection $\tau:V(G)\rightarrow V(G^*)$ preserving adjacency. The map $\tau$ is called an \emph{involution} if it satisfies the following:

1) $v\notin \tau (v)$ for every $v\in V$ and 

2) $u\in \tau (v)\iff v\in \tau (u)$ 

\begin{figure}
\centering
\includegraphics[width=.34\linewidth]{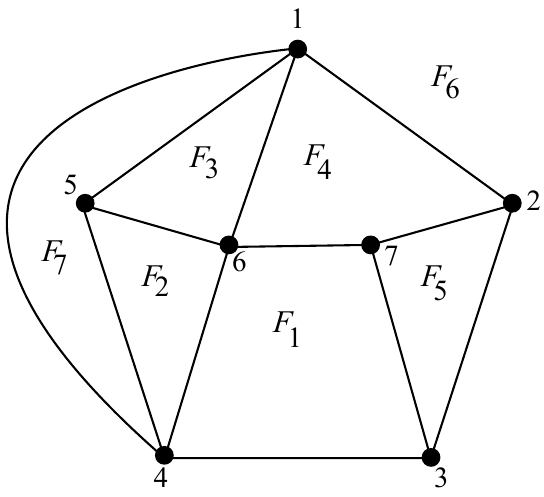}
\caption{Involutive graph with $\sigma(i)=i'$.}
\label{FigInv}
\end{figure}

A self-dual polyhedral graph $G$ admitting an involution is called an \emph{involutive} polyhedral graph (see \cite{montejano2020graphs}), see Figure \ref{FigInv}. 

Note that $\tau(v)$ can be thought as a face of $G$ (called \emph{dual} face of $v$, and denoted by $F_v$). It is easy to verify that for any edge $\{a, b\}\in E$, there is an other edge $\{x, y\}\in E$ such that $\tau(a)\cap\tau(b)=\{x, y\}$ and $\tau(x)\cap\tau(y)=ab$. We will write $\tau(\{a,b\})=\{x, y\}$ and call them \emph{dual edges}. Since the vertices of a Reuleaux polyhedron are in extremal configuration for the Vázsonyi problem, an involutive map exists and it is actually the canonical involution defined above in Theorem \ref{dualedges}. Hence, the 1-skeleton of a Reuleaux polyhedron is an involutive polyhedral graph.

Let $G=(V,E)$ be an involutive polyhedral graph and let $a,x\in V$. We say that $[a,x]$ is a \emph{diagonal} of $G$ if $x\in \tau (a)$. We define the \emph{diagonal} graph $\DDiag_G$ arising from $G$,  as the graph where the set of vertices is $V$ and set of edges is the set of all the diagonals of $G$. We observe that our diagonal graph correspond to the diameter graph used in \cite{montejano2020graphs}. We rather prefer to use the term diagonal to insist that it arises from the involutive map of the abstract graph.  In \cite{montejano2020graphs}, the authors studied involutive graphs from a more geometric point of view (in connection with \emph{metric mappings} and \emph{metric embeddings}) and thus the term diameter seems more appropriate. 

In \cite{montejano2020graphs},  the following was stated

\begin{conjecture}\label{ConjReul} \cite{montejano2020graphs}
Every involutive polyhedral graph $G=(V,E)$ is isomorphic to the 1-skeleton of a Reuleaux polyhedron $R(S)$ for some set of points $S$. 
\end{conjecture}

If this conjecture were true then we would have that $\DDiag_G$ is isomorphic to  $\DDiam_S$. Indeed, in such a case, there is a bijection $f:V \rightarrow S$ such that $[x,y]$ is a diagonal in $G$ if and only if the distance between $f(x)$ and $f(y)$ (vertices in the realization of $R(S)$) is equal to $\DDiam_S$.  Conjecture \ref{ConjReul} will be discussed further in the last section.

\begin{figure}
\centering
\includegraphics[width=.3\linewidth]{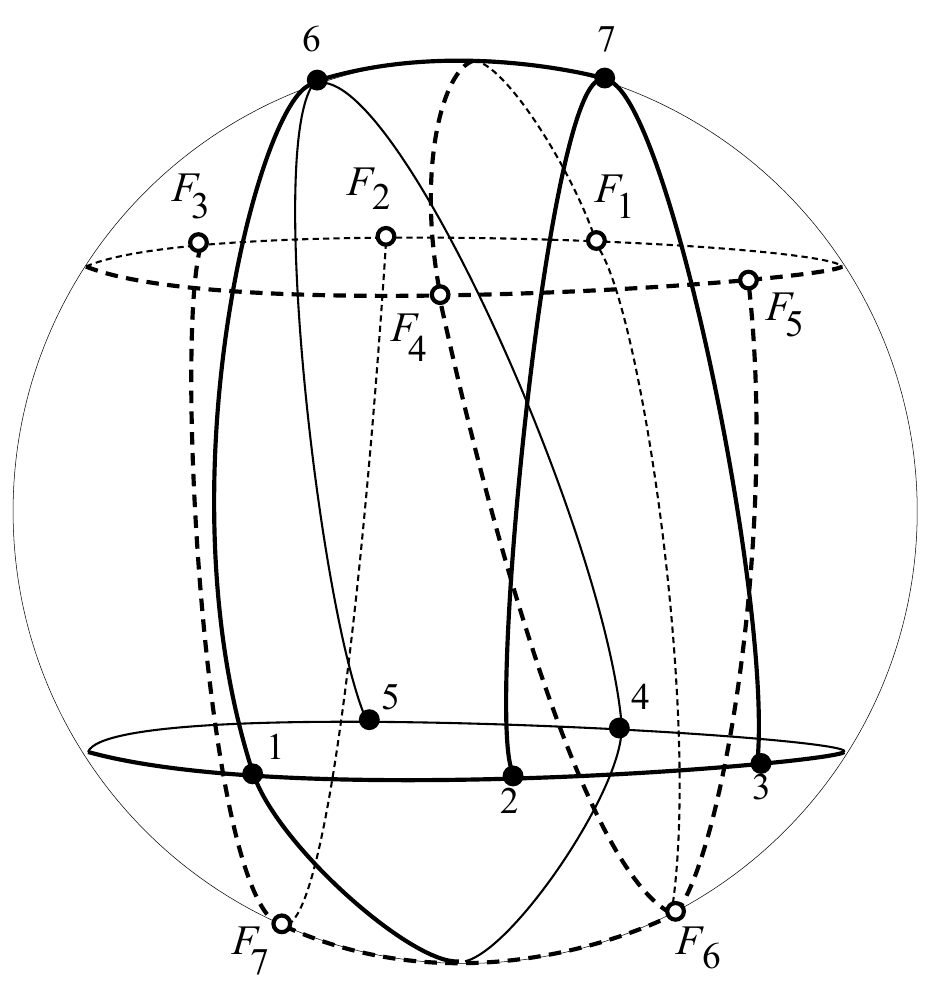}
\caption{Cellular embedding of graphs $H$ given in Figure \ref{FigInv}, (bold edges) and $H^*$ (dotted edges).}
\label{FigEmb}
\end{figure}

By Whitney's work \cite{whitney19332}, it is known that any polyhedral graph $G$ can be drawn in the plane or in the 2-sphere (in this case, $G$ is said to be a \emph{map}, that is, a graph cellularly embedded in $\mathbb{S}^2$) essentially in a unique way.  Montejano, Ram\'irez Alfons\'in and Rasskin \cite{montejano2022self} proved that any involutive polyhedra is \emph{antipodally self-dual}, that is, there are maps $\hat G$ and $\hat G^*$ that can simultaneously be cellularly embedded in $\mathbb{S}^2$ such that $-\hat{G}=\hat{G^*}$ where $-\hat G$ is the map consisting of the set of points $\{-x\in\mathbb{S}^2 \ |\  x\in\hat{G}\}$. In other words, if $x\in V(\hat{G})$ then $-x\in V(\hat{G}^*)$ and any edge $e^*\in E(\hat G^*)$ (respectively face $f^*\in F(\hat G^*)$) is antipodally embedded in $\mathbb{S}^2$ with respect to $e\in E(\hat G)$ (respectively $f\in F(\hat G)$), see Figure \ref{FigEmb}. 

\begin{figure}
\centering
\includegraphics[width=.3\linewidth]{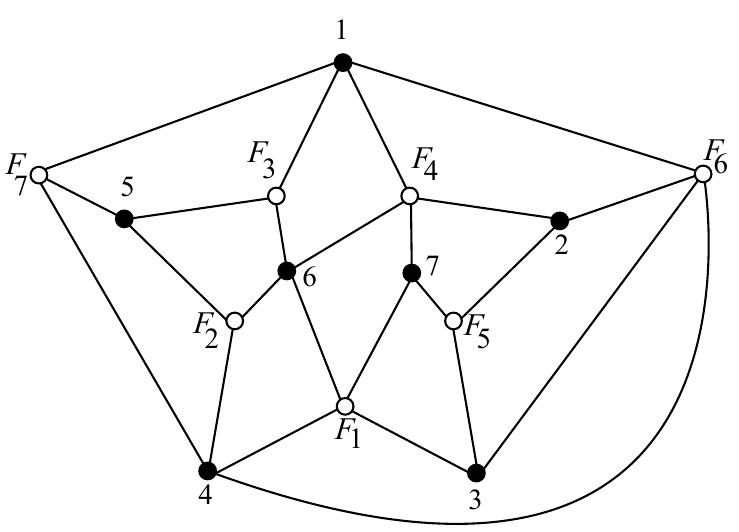}
\caption{Incident graph $I(H)$ with $H$ given in Figure \ref{FigInv}.}
\label{FigInc}
\end{figure}

Let $I(G)$ be the \emph{incidence} graph of the planar graph $G$. We recall that the set of vertices of $I(G)$ is given by $V(G)\cup V(G^*)$ and $\{v,w\}$ is an edge of $I(G)$ if $v\in V(G), w\in V(G^*)$ and $v\in F_{w}$ where $F_w$ is the face in $G$ corresponding to $w$, see Figure \ref{FigInc}. 

\begin{figure}
\centering
\includegraphics[width=.3\linewidth]{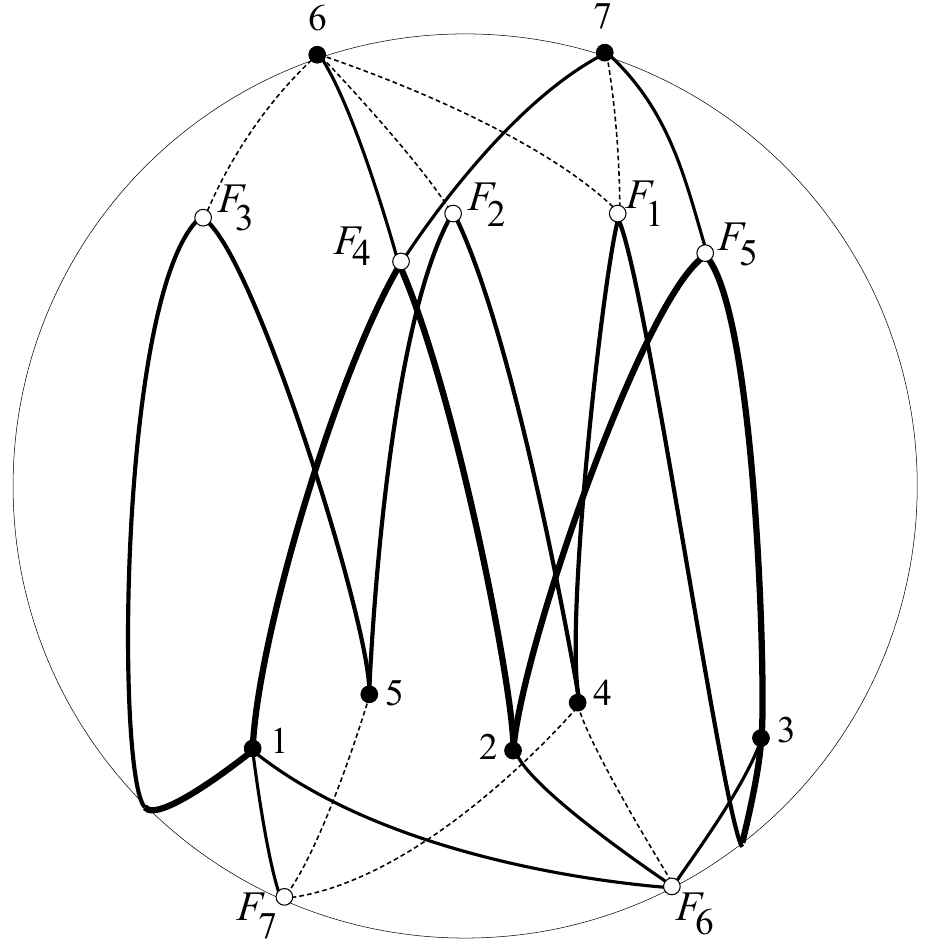}
\caption{Antipodal symmetric map $I(H)$ induced by the antipodal self-dual map given in Figure \ref{FigEmb}. A symmetric cycle is given by the bold edges}
\label{EmbInc}
\end{figure}

Notice that if $G$ is antipodally self-dual the embedding such that $-\hat{G}=\hat{G^*}$ naturally induce a {\em self-antipodally} map of $\hat I(G)$, that is, an embedding of $I(G)$ in $\mathbb{S}^2$ such that $\hat I(G)=-\hat I(G)$, that is, if $x\in V(\hat I)$ (respectively $x\in E(\hat I)$ or $x\in F(\hat I)$) then $-x\in V(\hat I)$ (respectively $-x\in E(\hat I)$ or $-x\in F(\hat I)$). We thus have that any vertex, edge or face in the antipodally symmetric map $\hat I$ is antipodally embedded in  $\mathbb{S}^2$ with respect to another vertex, edge and face of $\hat I$ (see \cite[Lemma 1]{montejano2022self}), see Figure \ref{EmbInc}. 

By a \emph{symmetric cycle} $C$ of $G$, we mean that there is an automorphism $\sigma$ of $G$ such that $\sigma(C) = C$ and $\sigma(int(C)) = ext(C)$, that is, $\sigma$ sends the cycle $C$ to itself and gives an isomorphism between the \emph{interior} of $C$ and the \emph{exterior} of $C$. In \cite[Lemma 1]{montejano2022self}, it was proved that if $G$ is an antipodally self-dual map then $I(G)$ is antipodally symmetric, see Figure \ref{EmbInc}.

\begin{figure}
\centering
\includegraphics[width=.3\linewidth]{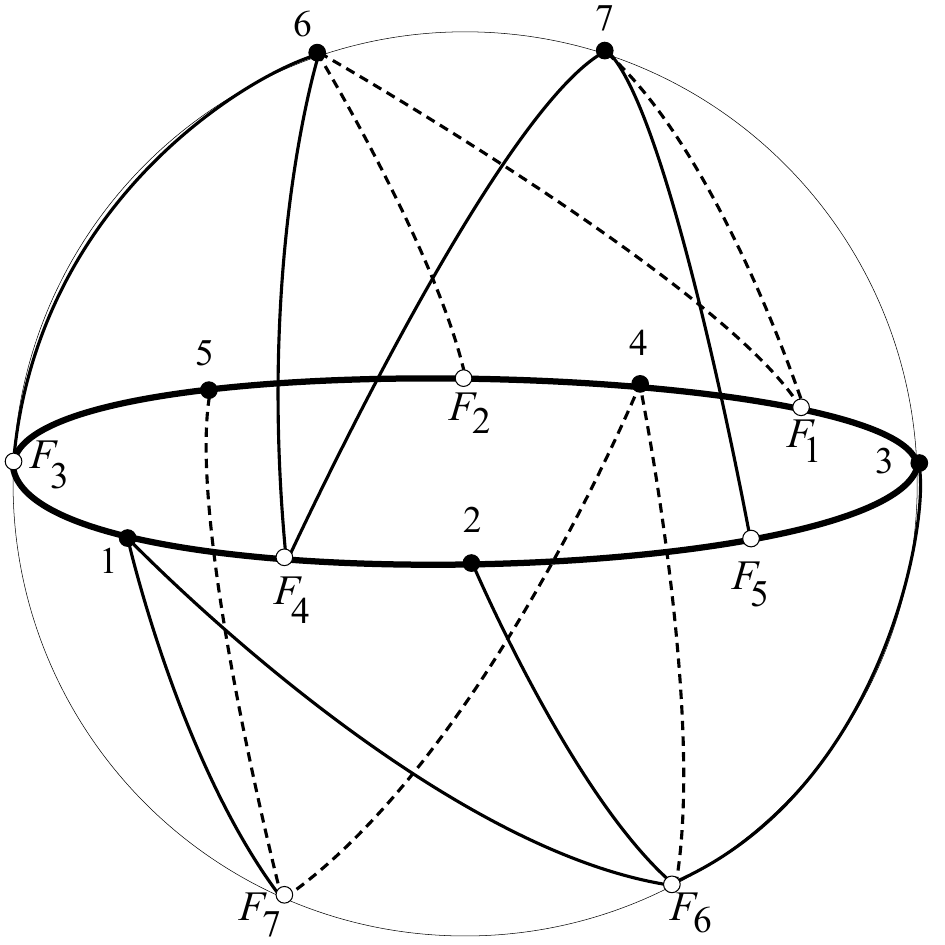}
\caption{Antipodal symmetric map $I(H)$ with a symmetric cycle (bold edges) in the equator with vertices forming a regular polygon.}
\label{EmbIncS}
\end{figure}

Furthermore, in \cite[Theorem 1]{montejano2022self} it was proved that if $G$ is an antipodally self-dual map then there is a \emph{symmetric cycle} $C_I$ with $2r$ vertices in $I(G)$ and $r$ odd. We shall denote by $\EEmbed(I(G))$ such embedding with $C_I$ placed along the equator of $\mathbb{S}^2$. This can be done keeping the antipodality of the embedding. Indeed, Once the cycle is on the equator we just draw $int(C)$ in the Northern hemisphere of $\mathbb{S}^2$ and then we draw $ext(C)$ in the Southern hemisphere in an antipodally fashion, that is, for each vertex (or edge) in $int(C)$ there is its antipodal vertex (or edge) in $ext(C)$. For more details about the properties of this embedding we refer the reader to \cite{montejano2022self}, see Figure \ref{EmbIncS}. 

 The notion of symmetric cycle in maps has already been used in other contexts, for instance, to study knot theory problems [\cite{montejano2023self2}, \cite{montejano2022selfdual3}].

\section{The Key Lemma}\label{KeyLemma}

This section is devoted to proving the following lemma that plays a central role throughout this paper.

\begin{lemma}\label{critical}
 Let $G$ be an involutive polyhedral graph. Then, $\DDiag_G$ is 4-critical, that is, it is vertex 4-chromatic and the removal of any vertex decreases its chromatic number.    
 \end{lemma}

In order to prove the above lemma, we first establish a number of important properties needed as basic bricks for its proof.
 
Let $G$ be an involutive graph. We shall consider the above aforementioned antipodal embedding $\EEmbed(I(G))$ in $\mathbb{S}^2$ where the symmetric cycle $C_I$ is minimal, that is, with a minimal number of edges. We suppose that $int(G)$ and $ext(G)$ are drawn in the open Northern and the Southern hemispheres (denoted by 
$\mathbb{S}^2_{\text{\tiny N}}$ and $\mathbb{S}^2_{\text{\tiny S}}$) respectively.

\begin{figure}
    \centering
    \begin{subfigure}[b]{0.3\textwidth}
        \includegraphics[width=\textwidth]{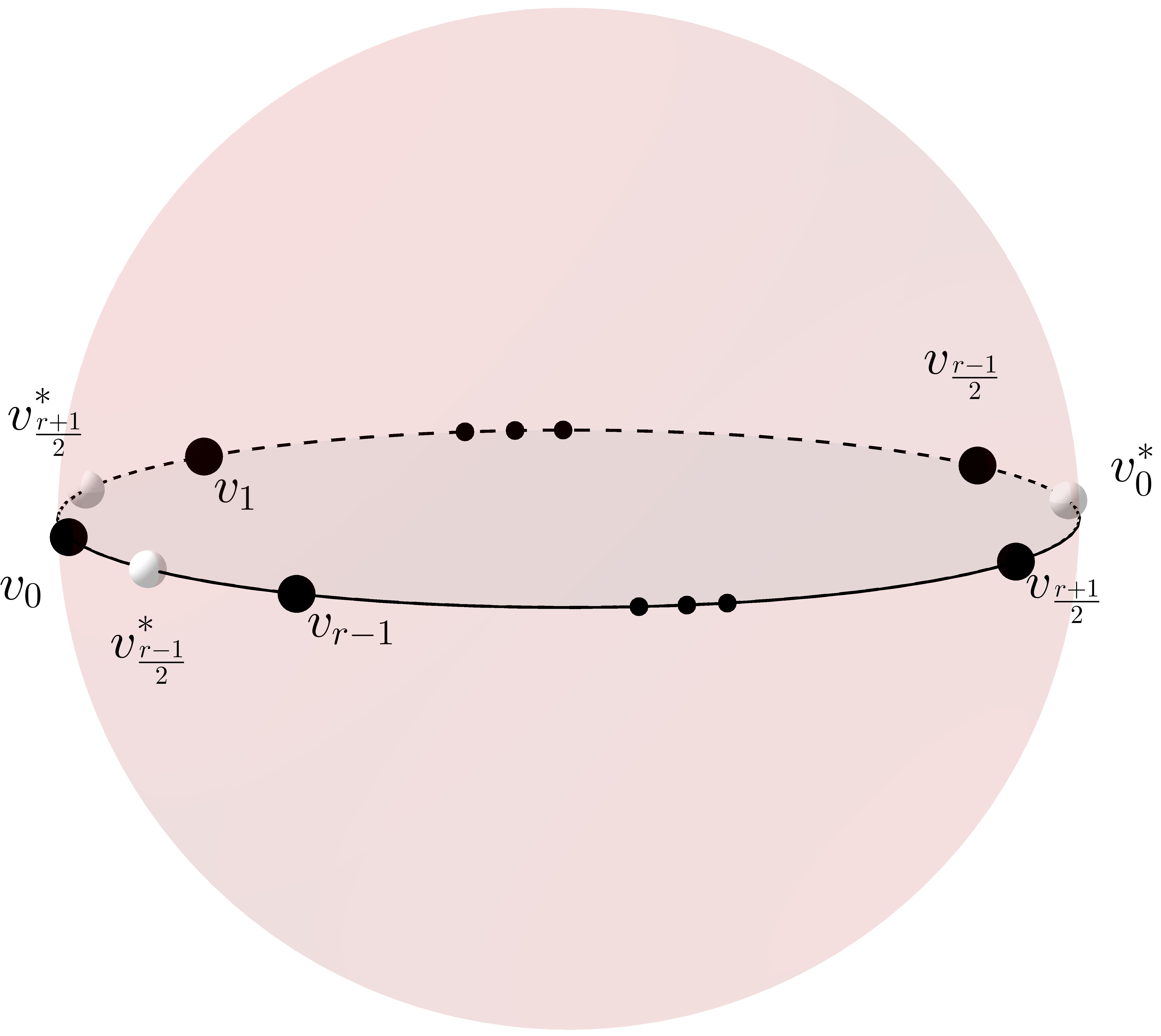}
        \caption{Cycle $C_I$ with $r=7$.}
        \label{CycleCI}
    \end{subfigure}
     \hfill
     \begin{subfigure}[b]{0.3\textwidth}
        \includegraphics[width=\textwidth]{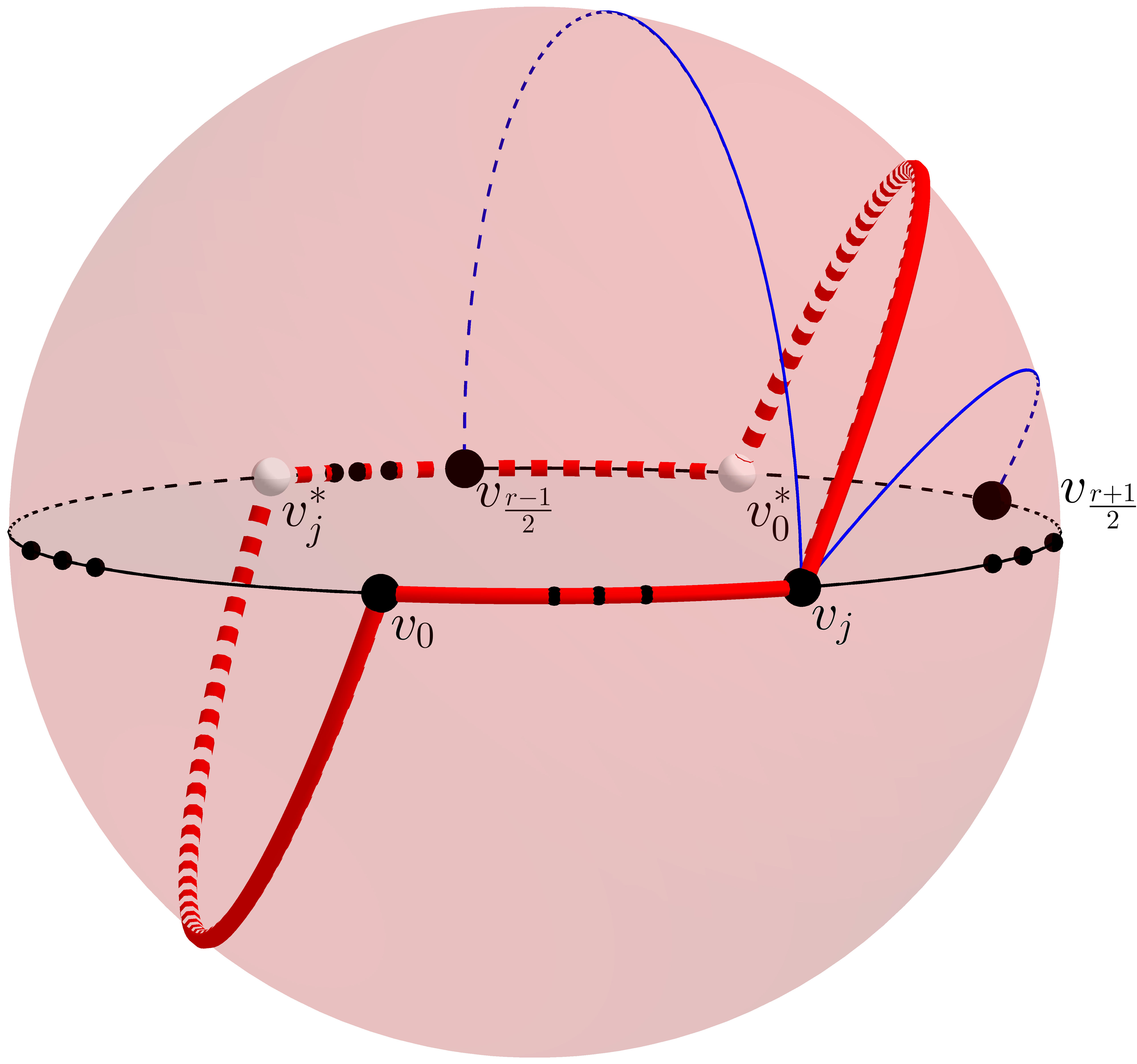}
        \caption{Red edges inducing the new shorter symmetric cycle $C_I$. Blue arcs are edges in $G$.}
        \label{ExtraFace}
    \end{subfigure}
     \hfill
    \begin{subfigure}[b]{0.3\textwidth}
        \includegraphics[width=\textwidth]{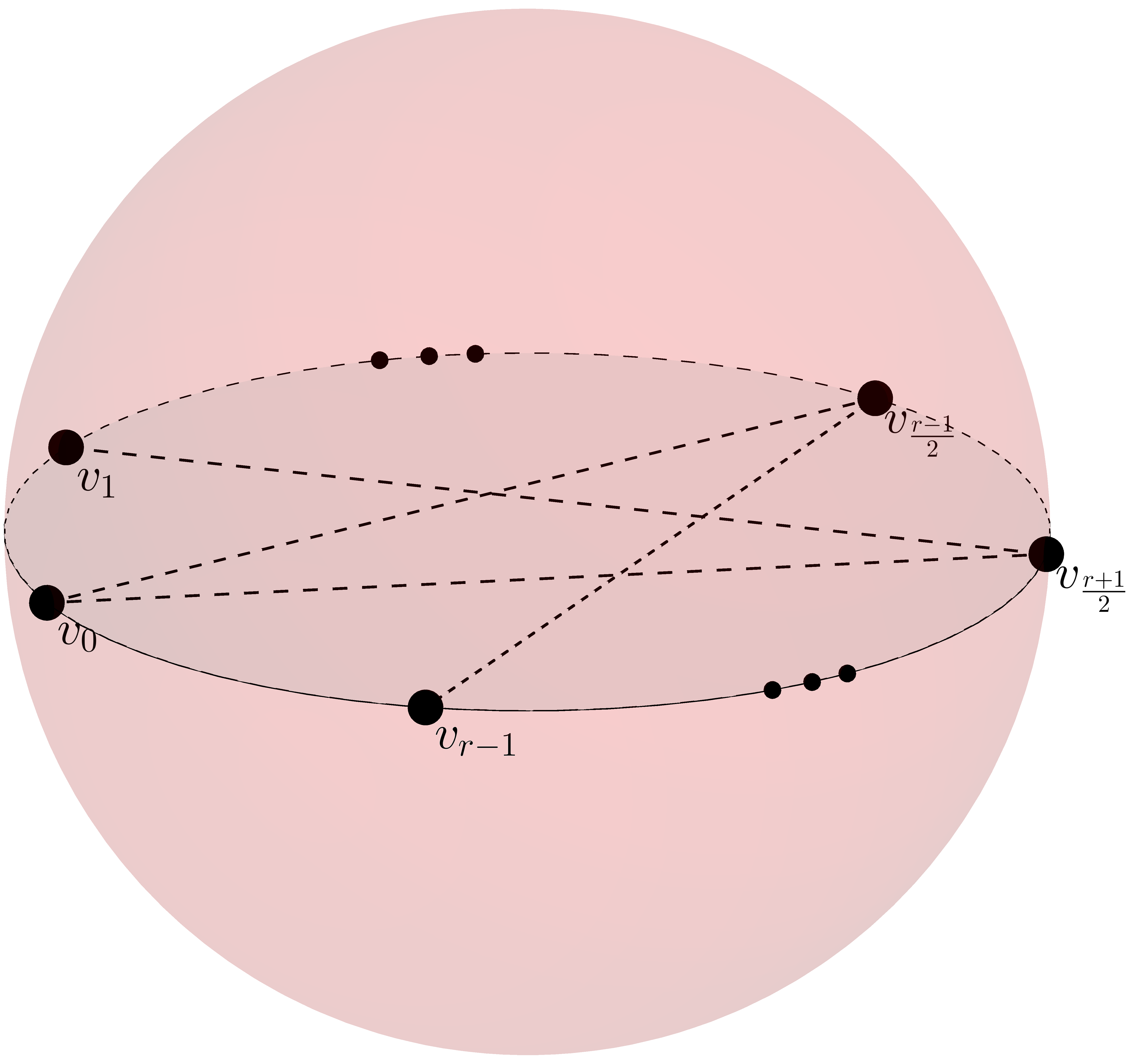}
        \caption{Cycle $C_D$.}
        \label{CycleCD}
    \end{subfigure}
    \caption{Edges of $I(G)$ in black and edges of $G$ in blue}
    \label{Properties}
   
\end{figure}

{\bf [P1]} We suppose that $|C_I|=2r$ where $r$ is an odd integer. We label with color black (respectively,  white) the vertices $v_0,\dots ,v_{r-1}$ (respectively $v_0^*,\dots ,v_{r-1}^*$) clockwise around the equator which are the vertices in $C_I\cap G$ (respectively in $C_I\cap G^*$).  Since vertex $v_i$ is antipodally embedded to $v^*_i$ so $C_I$ is cyclically labeled as follows $\{v_0,  v_{\frac{r+1}{2}}^*, v_1, \dots, v_{\frac{r-1}{2}}, v_0^*, v_{\frac{r+1}{2}}, v_1^* \dots, v_{\frac{r-1}{2}}^* \}$, see Figure \ref{CycleCI} 

{\bf [P2]} We claim that any $v_i$ is adjacent to exactly two vertices of $C_I$ in $\DDiag_G$. We may show this for $v_0$ (the argument is the same for any $v_i$). Recall that $v_0^*$ is a vertex of $G^*$ representing the face $F_{v_0}$ determined by the duality isomorphism, say $\sigma$, of $G$ that sends $v_0$ to $v_0^*$. Since $v_{\frac{r-1}{2}}$ and $v_{\frac{r+1}{2}}$ are adjacent to  $v_0^*$ (in the cycle $C_I$) then, by definition of $I(G)$, these vertices belong to the face $F_{v_0}$. Since $\sigma$ is involutive then all the vertices of $F_{v_0}$ are adjacent to $v_0$ in $\DDiag_G$.

Now, suppose that there is another $v_j$, $j\neq \frac{r-1}{2}, \frac{r+1}{2}$ adjacent to $v_0$. The latter means that $v_j$ is also in the face $F_{v_0}$ and therefore there must also exists an edge joining $v_j$ and $v_0^*$ in $G_I$, see Figure \ref{ExtraFace}.

Since $I(G)$ is antipodally symmetric then, there is also an edge joining $v_j^*$ and $v_0$. We way construct the cycle $C'_I=v_0,v_j^*, [v_j^*;v_0^*],v_0^*,v_j, [v_j;v_0]$ where $[a;b]$ denotes the path along the equator joining $a$ and $b$ without intersecting any other previous vertex in $C'_I$.
By the antipodality of $I(G)$, we have that $C'_I$ induce a symmetric cycle of $I(G)$ with $|C_I'|<|C_I|$, which is a contradiction to the minimality of $C_I$, see Figure \ref{ExtraFace}.

{\bf [P3]} By [P2],  the degree of each vertex $v_i$ of $C_I$ in $\DDiag_G$ is equals two. In other words, $v_i$ form two diagonals with the two vertices adjacent to $v_i^*$ in $C_I$. Since $r$ is odd then the set of all these couple of diagonals  form a cycle $C_D$ in $\DDiag_G$. $C_D$ is a star with $r$ vertices in $C_I$. For commodity, we preserve the same vertex labels of $C_I$ , given by the order of appearance around the equator for $C_D$, see Figure \ref{CycleCD}

\begin{figure}
    \centering
    \begin{subfigure}[b]{0.3\textwidth}
        \includegraphics[width=\textwidth]{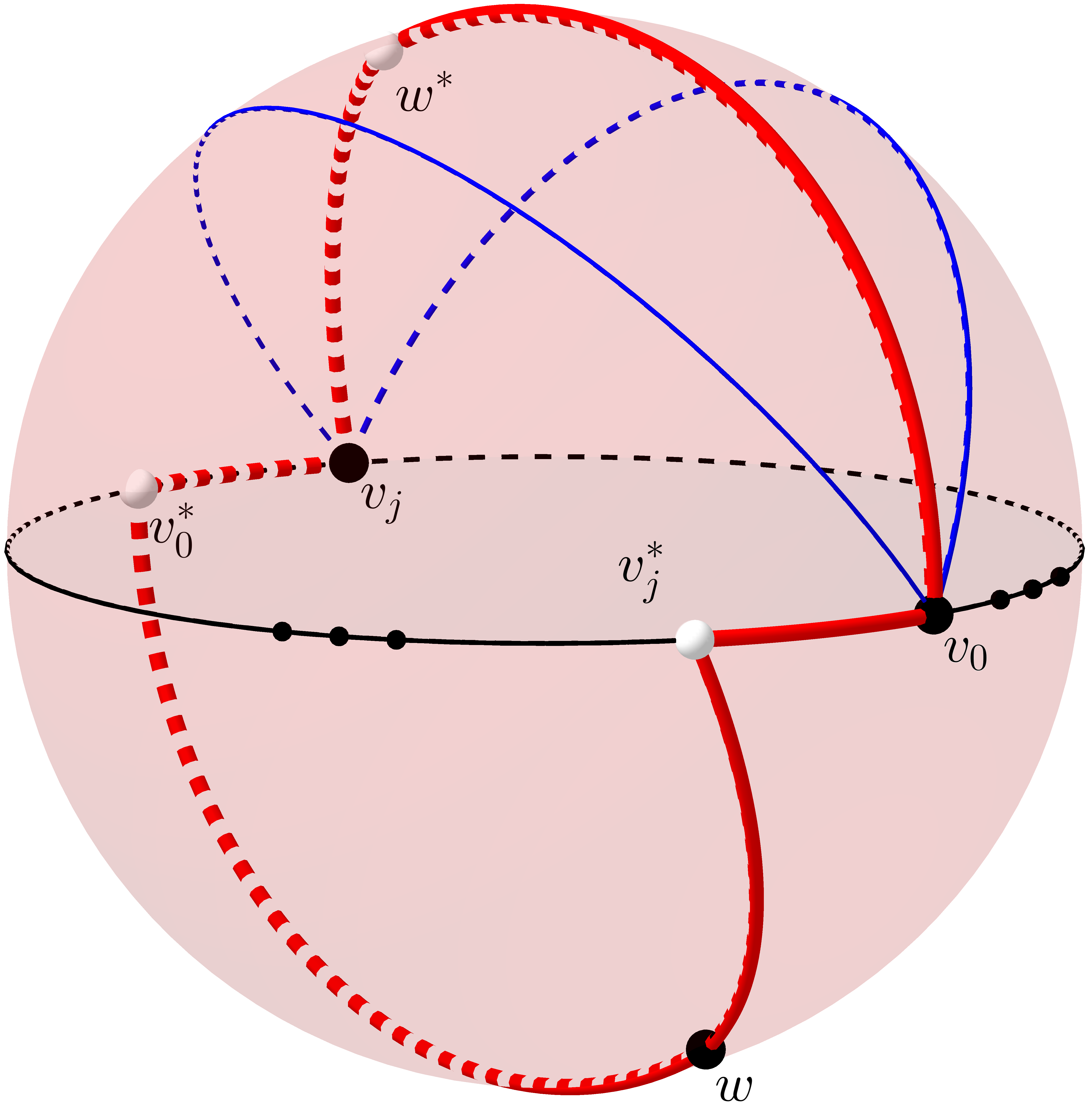}
        \caption{Red edges inducing the new shorter symetric cycle $C_I'$.}
        \label{NotConsecutive}
    \end{subfigure}
     \hfill
    \begin{subfigure}[b]{0.3\textwidth}
        \includegraphics[width=\textwidth]{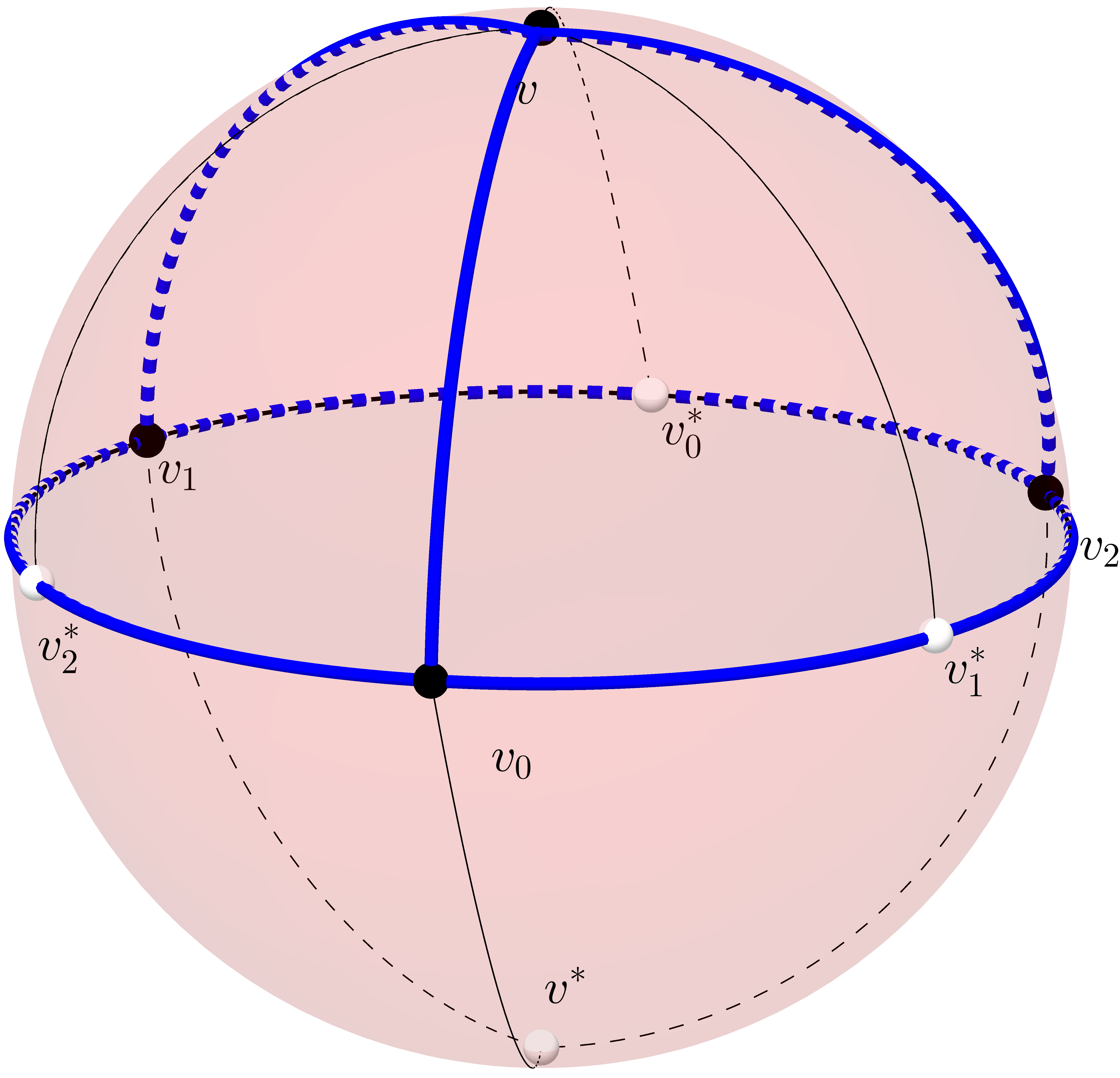}
        \caption{The Tetrahedron.}
        \label{Tetrahedron}
    \end{subfigure}
     \hfill
    \begin{subfigure}[b]{0.3\textwidth}
        \includegraphics[width=\textwidth]{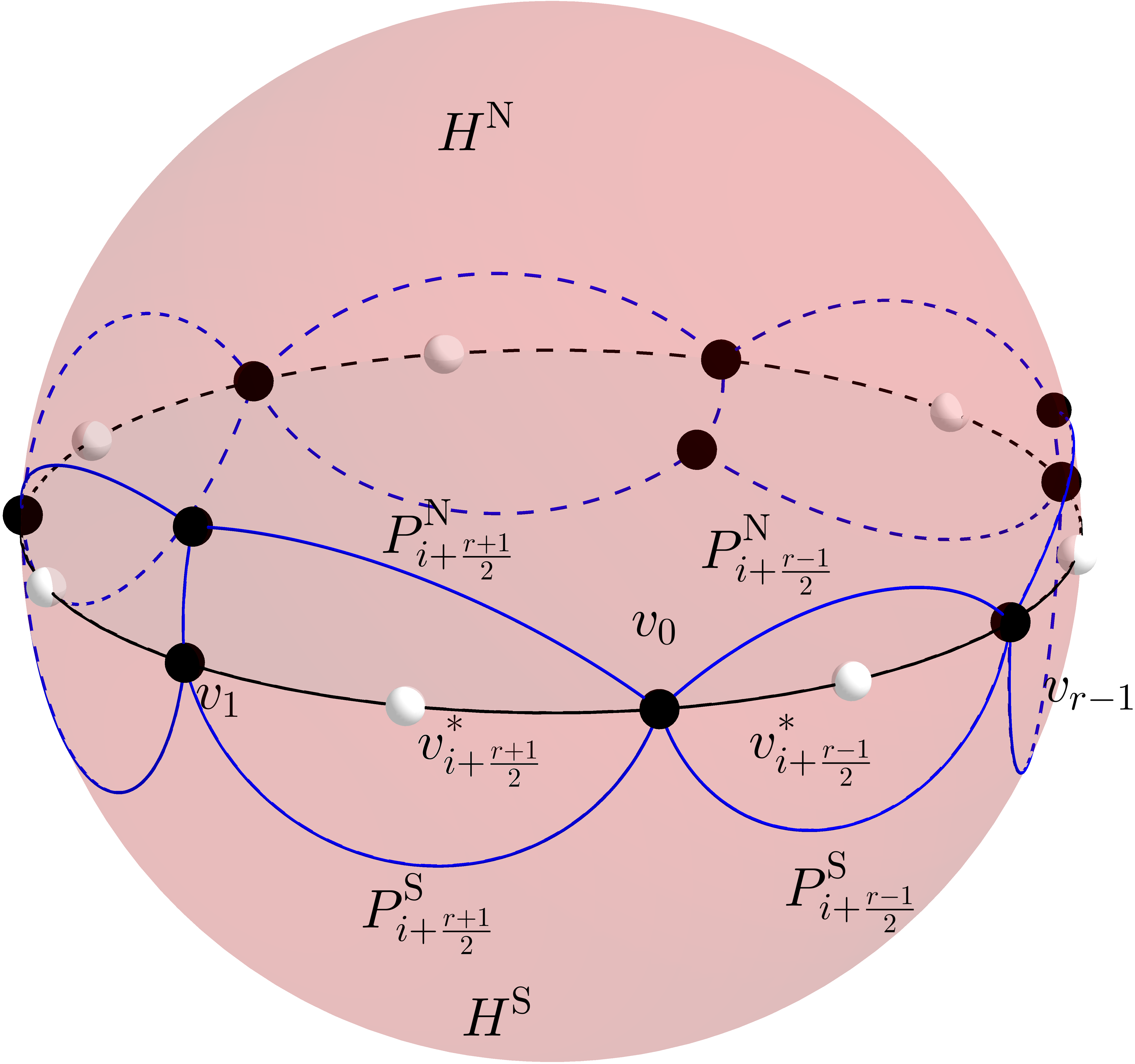}
        \caption{$H^{\text{N}}, H^{\text{S}}, P^N_i$ and $P^S_i$}
        \label{PHCircuit}
    \end{subfigure}
    \caption{Edges of $I(G)$ in black and edges of $G$ in blue}
    \label{Properties2}
   
\end{figure}

{\bf [P4]} We claim that there is not face of $G$ containing two non-consecutive vertices of $C_D$ (recall that consecutive is with respect to the order of appearance around the equator and not in the order of appearance while traveling through $C_D$).  We proceed by contradiction, suppose that there is a face $F_w$ containing two 
non-consecutive vertices, say $v_0$ and $v_j$.  We thus have that the vertex $w^*$, representing the dual face $F_w$ in $I(G)$, must be adjacent to both $v_0$ and $v_j$. By antipodality, we also have that $w$ is adjacent to both $v_0^*$ and $v_j^*$. We may thus construct a symmetric cycle $C'=[v_0;v_j^*],v_j^*,w,v_0^*,[v_0^*;v_j],v_j, w^*,v_0$ with $|C'|<|C_I|$, which is a contradiction to the minimality of $C_I$, see Figure \ref{NotConsecutive}.

{\bf [P5]} Notice that a face $F$ of $G$ can never contain four or more vertices of $C_D$, otherwise $F$ would have at least two non-consecutive vertices of $C_D$ which, by [P4],  it is impossible.

There might exist a face $F$ containing exactly three consecutive vertices of $C_D$, in this case, $G$ is actually the tetrahedron. Indeed, since the vertices are consecutive then  $C_D$ consists of three vertices and thus the drawing of $I(G)$ consist of six vertices in the equator (three black and three white appearing alternating) with one black vertex say in $int(G)$ joined to the three white  vertices in the equator and one white vertex in $ext(G)$ (representing the face $F$) joined to the three black vertices in the equator. We thus have that $G$ consist of 4 black vertices forming a tetrahedron, see Figure \ref{Tetrahedron}

{\bf [P6]} Let us consider the embedding of $G$ in $\mathbb{S}^2$, say $\EEmbed(G)$,  induced by the embedding of $I(G)$.  By the symmetry of $C_I$,  the only faces in  $\EEmbed(G)$ lying in both $\mathbb{S}^2_{\text{\tiny N}}$ and $\mathbb{S}^2_{\text{\tiny S}}$ at once  are the faces corresponding to each black vertex in $C_I$. Any other face completely lies in either of the hemispheres, see Figure \ref{PHCircuit}.

{\bf [P7]} Recall that $F_{v_i}$ is the dual face of $v_i$ represented by vertices $v_i^*$. We define $P_i^{\text{\tiny N}}$ (respectively $P_i^{\text{\tiny S}}$) as the path going from $v_{i+\frac{r-1}{2}}$  to $v_{i+\frac{r+1}{2}}$ for  each $i=0,\dots, \frac{r+1}{2}$ (sum $\bmod\ r$) through the vertices of $F_{v_i}$ appearing in $\mathbb{S}^2_{\text{\tiny N}}$ (respectively in $\mathbb{S}^2_{\text{\tiny S}}$). 

We also let $H^{\text{N}}$ (respectively $H^{\text{S}}$) be the union of all $P_i^{\text{\tiny N}}$ (respectively all $P_i^{\text{\tiny S}}$), see Figure \ref{PHCircuit}.

{\bf [P8]} Since $G$ is a polyhedral graph (and thus simple) then any pair of faces shares at most one edge. Therefore, we may have repeated consecutive edges in $H^{\text{N}}$ (or $H^{\text{S}}$) if $F_{v_i}$ and $F_{v_{i+1}}$ share an edge, see Figure \ref{PHCircuit}.

{\bf [P9]} Notice that $H^{\text{N}}$ (respectively $H^{\text{S}}$) induces a path of $G$ separating all the faces completely contained in $\mathbb{S}^2_{\text{\tiny N}}$ (respectively in $\mathbb{S}^2_{\text{\tiny S}}$) from the rest of the faces, see Figure \ref{PHCircuit}.

We may now prove Lemma 1.

\begin{Proof1}
    
By \cite[Theorem 3]{montejano2020graphs}, $\chi(\DDiag_G)=4$. We shall show that $\chi(\DDiag_G\setminus \{v\})=3$ for any $v\in V(\DDiag_G)$. In order to prove this,  for each $v\in V(G)$ we will show that there is always a map $c: V(G)\rightarrow \{0,1,2,3\}$ from the vertices of $G$ to colors 0,1,2 and 3 inducing a proper coloring with $c(v)=3$ and $c(v)\neq c(u)$ for all $u\neq v$.

We have that either $v$ is a vertex in $V(C_D)$ or it lies in an hemisphere. Let us see each of these two cases.

{\bf Case 1)} Let $v\in V(C_D)$. Without loss of generality, we may take $v=v_0$ (in the labeling of $C_D$). We have that the dual face $F_{v_0}$ contains at least three vertices, say $v_{\frac{r-1}{2}}, v_{\frac{r+1}{2}}$ (see {\bf [P7]}) and $u$. Without loss of generality,  we may assume that $u$ lies in $\mathbb{S}^2_{\text{\tiny N}}$. 

Let us observe that, by definition of the paths $P_{i}^{\text{\tiny N}}$  (see {\bf [P7]}), $u\in P_{\frac{r-1}{2}}^{\text{\tiny N}}$. We will use this fact later on in the Subcase 1.2 below.

Let $A[v_0,v_{\frac{r+1}{2}}]$ (respectively $A[v_0,v_{\frac{r-1}{2}}]$) be the vertices in the arc of the equator between $v_0$ and $v_{\frac{r+1}{2}}$ not containing $v_{\frac{r-1}{2}}$ (respectively the arc between $v_0$ and $v_{\frac{r-1}{2}}$ not containing $v_{\frac{r+1}{2}}$), see Figure \ref{Case1Coloration}.

We color the vertices of $G$ as follows. 

\begin{itemize}
    \item $c(v=v_0)=3$,
    \item $c(x)=2$ if  $x\in A[v_0,v_{\frac{r+1}{2}}]\backslash\{v_0\}$,
    \item $c(x)=1$ if  $x\in A[v_0,v_{\frac{r-1}{2}}]\backslash\{v_0\}$,
    \item $c(x)=0$  if $x$ lies in $\mathbb{S}^2_{\text{\tiny N}}$,
    
\end{itemize}
see Figure \ref{Case1Coloration}

\begin{figure}
    \centering
    \includegraphics[scale=0.3]{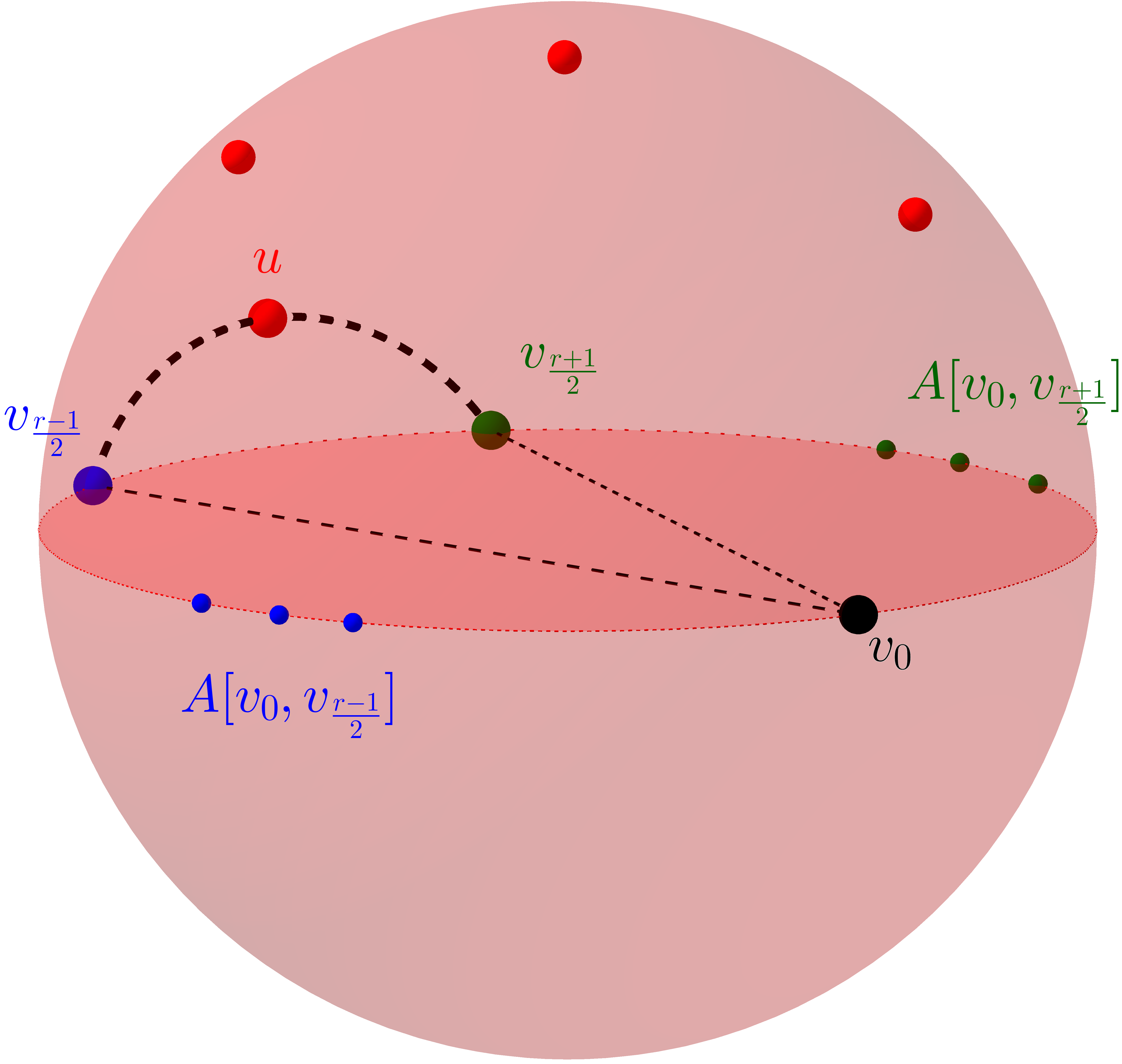}
    \caption{The blue vertices are color 1, the green ones are color 2 and the red ones are color 0.}
    \label{Case1Coloration}
   
\end{figure}
 
We first observe that the vertices of an edge in $C_D$ have different colors (and thus colored properly). Moreover, since there is no edge of $\DDiag_G$ between two vertices in $\mathbb{S}^2_{\text{\tiny N}}$ (all the neighbors of the vertices in $\mathbb{S}^2_{\text{\tiny N}}$ in $\DDiag_G$ lie in $\mathbb{S}^2_{\text{\tiny S}}$ our coloring works so far.

We finally need to color each vertex lying in $\mathbb{S}^2_{\text{\tiny S}}$. Let $w$ be a vertex of $G$ in $\mathbb{S}^2_{\text{\tiny S}}$ and let $F_w$ be its dual face lying in $\mathbb{S}^2_{\text{\tiny N}}$. We claim that at most two out of the three colors 0,1 and 2 can be used for the vertices in $F_w$.  If this is the case, we may then color vertex $w$ with a color different from  0,1 and 2.
We prove the claim by contradiction. Let us suppose that the three colors 0,1 and 2 are used in the vertices of $F_w$. 
If colors 1 and 2 are used then $F_w$ must have two vertices of $C_D$. By {\bf [P4]}, these vertices cannot be non-consecutive, and therefore the only choice for these vertices to be in $F_w$ are $v_{i+\frac{r-1}{2}}$ and $v_{i+\frac{r+1}{2}}$. 

We have two subcases.

{\bf Subcase 1.1)} We suppose that  $u\in F_w$. We claim that $F_w=F_{v_0}$. Indeed, $F_w$ and $F_{v_0}$ have three common vertices and since any two faces share at most one edge (since $G$ is 3-connected) then the only way for this to happen is if $F_w=F_{v_0}$. However, the latter implies that $w=v_0$, contradicting the fact that $w$ is in $\mathbb{S}^2_{\text{\tiny S}}$.

{\bf Subcase 1.2)} We suppose that $u\not\in F_w$. Since both faces $F_{v_0}$ and $F_w$ passe through $v_{\frac{r-1}{2}}$ and $v_{\frac{r+1}{2}}$ then $F_w$ must contain $F_{v_0}$, in particular, $F_w$ contains $P_{\frac{r-1}{2}}^{\text{\tiny N}}$, see Figure \ref{Case1Coloration}.  As observed above, $u\in P_{\frac{r-1}{2}}^{\text{\tiny N}}$. We clearly have that any path connecting $u$ with any other vertex in the exterior of $F_w$ must go through either  $v_{\frac{r-1}{2}}$ or $v_{\frac{r+1}{2}}$, implying that these are cut vertices, contradicting the 3-connectivity of $G$.

{\bf Case 2)} Let $v$ be a vertex lying in $\mathbb{S}^2_{\text{\tiny N}}$ (the case when $v$ lies in $\mathbb{S}^2_{\text{\tiny S}}$ is analogous). We will first construct three vertex-disjoint paths joining $v$ with three different vertices of $C_D$. 

Let $w$ be a vertex in $\mathbb{S}^2_{\text{\tiny S}}$ (this vertex exists, otherwise $G$ would be the tetrahedron which is clearly 4-critical).
Since $G$ is 3-connected then, by Menger's theorem, there exist three vertex-disjoint paths $Q_0, Q_1$ and $Q_2$ joining $u$ to $w$. We clearly have that each of these paths must intersect $H^{\text{N}}$.  Let $h_i$ be the first vertex of $H^{\text{N}}$ hit by $Q_i$ from $v$ to $w$,  for each $i=0,1,2$. Suppose that $h_i$ is in  one of the $P_{v_i}^{\text{\tiny N}}$, we denote it by $P(h_i)$ for short. We observe that there are two ways to reach $C_D$ from $h_i$ : either by following the vertices of $P(h_i)$ appearing to the right of $h_i$ (denoted by $R_i$) or by following the vertices of $P(h_i)$ appearing to the left of $h_i$ (denoted by $L_i$). Notice that $R_i$ or $L_i$ may consist of only the vertex $h_i$, which is already a vertex in $C_D$.

These are the desired paths, (we write $Q_i[v,  h]$ the sub-path of $Q_i$ from $v$ to $h$) :
\begin{itemize}
    \item $\bar Q_0:=Q_0[v,h_0]\cup T_0$, where $T_0$ is either $R_0$ or $L_0$.
    \item $\bar Q_1:=Q_1[v,h_1]\cup T_1$, where $T_1$ is either $R_1$ or $L_1$. Notice that if $P(h_0)=P(h_1)$ then we can always take $T_1$ as the side not used in $T_0$.
    \item $\bar Q_2:=Q_2[v,h_2]\cup T_2$, where $T_2$ is either $R_2$ or $L_2$.
Notice that  if $P(h_0)=P(h_1)$ then $P(h_2)\neq P(h_0), P(h_1)$ otherwise there will be two $Q_i$'s with a common vertex (which is not possible since they are vertex-disjoint), see Figure \ref{DisjointPathsImposible}
\end{itemize}

 \begin{figure}
    \centering
    \begin{subfigure}[b]{0.4\textwidth}
        \includegraphics[width=\textwidth]{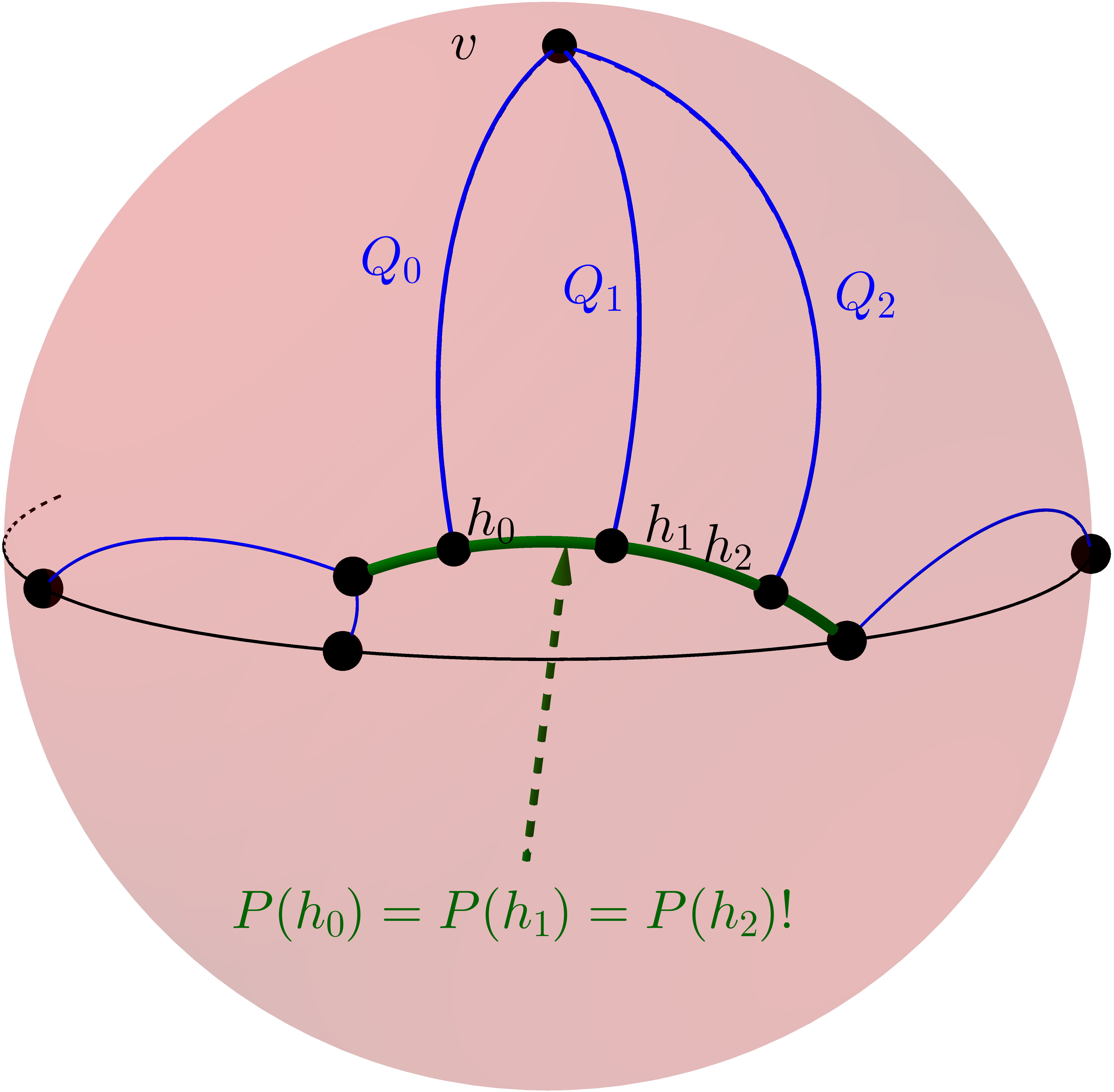}
        \caption{$Q_1$ would intersect either $Q_0$ or $Q_2$.}
        \label{DisjointPathsImposible}
    \end{subfigure}
     \hfill
    \begin{subfigure}[b]{0.4\textwidth}
        \includegraphics[width=\textwidth]{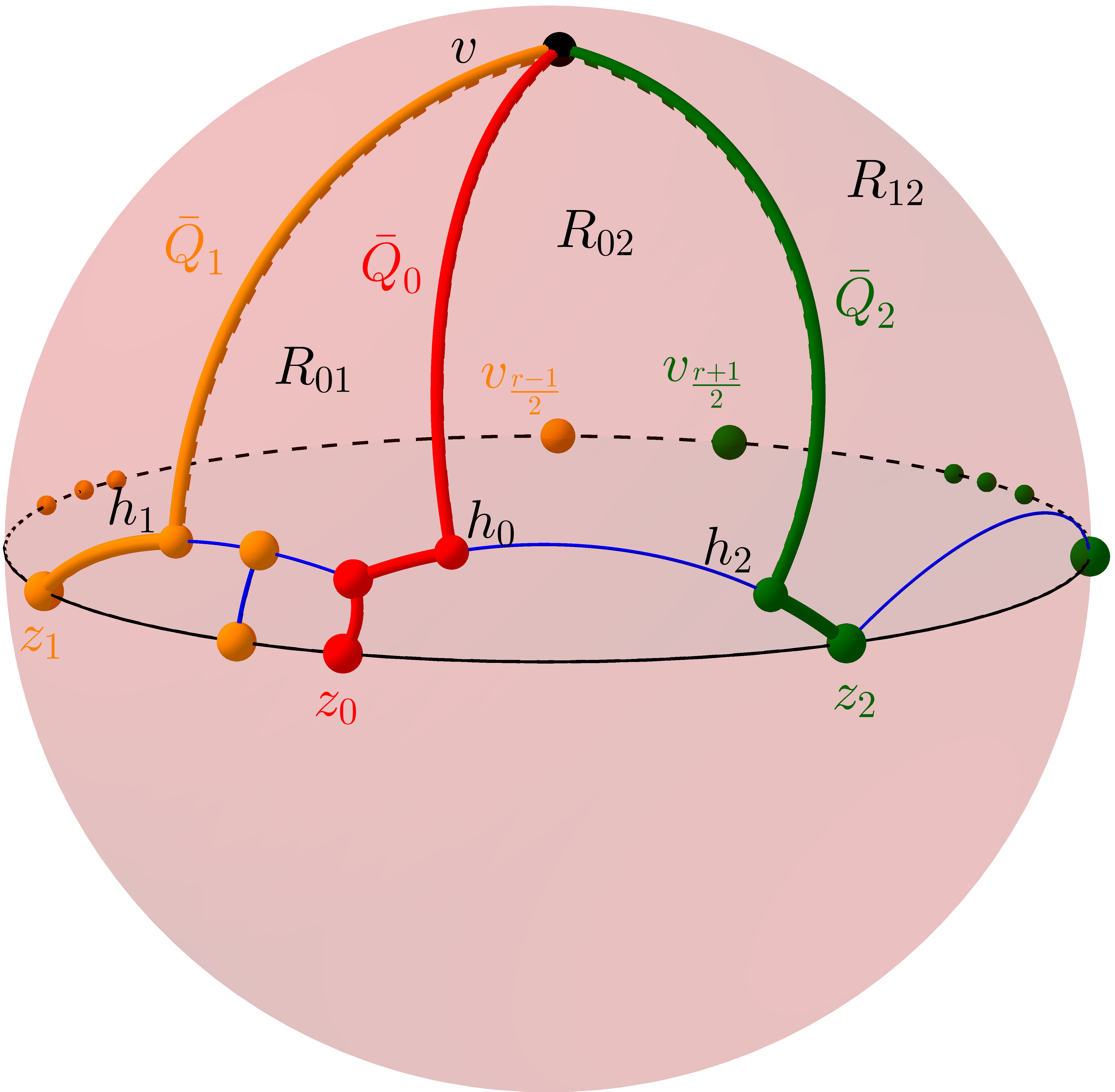}
        \caption{Coloration by $\bar Q_0, \bar Q_1$ and $\bar Q_2$.}
        \label{DisjointPaths}
    \end{subfigure}
    
    \caption{Division of the north hemisphere by the paths $\bar Q_0, \bar Q_1$ and $\bar Q_2$.}
    \label{PathsQ}
   
\end{figure}

Suppose that the vertices $v_i$ are placed in a $r$-regular polygon with $r$-odd,  all on the equator. Notice that since the polygon is regular,  then all its vertices lie on a common circle and since $r$ is odd then for each vertex,  there is a  unique diametrically opposite side. Let $z_i$ be the common vertex of $\bar Q_i$ and $C_D$. Draw a line $\ell$ (lying on the plane containing the equator) going through $z_0$ and perpendicular to the opposite side of the regular polygon. We may suppose that we have the situation in which $z_1$ and $z_2$ are in opposite sides of $\ell$. Otherwise, if both $z_1$ and $z_2$ are on the same side of $\ell$ then either $z_1$ is between $z_0$ and $z_2$ or $z_2$ is between $z_0$ and $z_1$.
If $z_1$ is between $z_0$ and $z_2$,  then we clearly have that the line $\ell'$ going through $z_1$ perpendicular to the opposite side in the regular polygon will leave $z_0$ and $z_2$ on different sides (similarly if $z_2$ were the middle vertex). 

Without loss of generality,  we may assume that $z_0=v_0$. Let $A[v_0,v_{\frac{r+1}{2}}]$ (respectively $A[v_0,v_{\frac{r-1}{2}}]$) be the vertices in the arc of the equator between $v_0$ and $v_{\frac{r+1}{2}}$ containing $z_1=z_x$ (respectively between $v_0$ to $v_{\frac{r-1}{2}}$ containing $z_2=z_y$), see Figure \ref{DisjointPaths}.

We begin coloring some vertices lying in $C_D$ and $\mathbb{S}^2_{\text{\tiny N}}$ as follows:

\begin{itemize}
    \item $c(v)=3$, 
    \item $c(v_0)=0$, 
    \item $c(x)=1$ for all vertex $x\in A[v_0,v_{\frac{r+1}{2}}]\setminus\{v_0\}$, 
    \item $c(x)=2$ for all vertex $x\in A[v_0,v_{\frac{r-1}{2}}]\setminus\{v_0\}$,
    \item $c(x)=0$ for all vertex $x\in \bar Q_0\setminus \{v\}$,
    \item $c(x)=1$ for all vertex $x\in \bar Q_1\setminus \{v\}$ and
    \item $c(x)=2$ for all vertex $x\in \bar Q_2\setminus \{v\}$, see Figure \ref{DisjointPaths}
\end{itemize}

Let us verify that this partial coloring is correct so far. We first remark that any vertex in $A[v_0,v_{\frac{r+1}{2}}]\setminus\{v_0\}$ (with color 1) is correctly colored since its neighbors are two opposite vertices lying in $A[v_0,v_{\frac{r-1}{2}}]$ having color 2 (similarly, for the vertices in $A[v_0,v_{\frac{r-1}{2}}]\backslash \{v_0\}$).

Let us check that the vertices in $\bar Q_i=Q_i[v,h_i]\cup T_i$ are all well colored. We notice that there is not problem with the colors of vertices in $Q_i[v,h_i]$ since all their neighbors (in $\DDiag_G$) are vertices in $\mathbb{S}^2_{\text{\tiny S}}$ (which are not colored yet). Let us now check the vertices of $T_i$. We will do so for $T_0$ (analogous arguments can be used to check that the vertices in both $T_1$ and $T_2$ are also properly colored).

We have that the vertices of $T_0=[h_0,  \dots ,v_0]$ (colored with color 0 since they are contained in $\bar Q_0$) is a subset of $P_{v_{\frac{r-1}{2}}}$,  which, in turn, as pointed out in {\bf [P7]}, is a subset of the dual face $F_{v_{\frac{r-1}{2}}}$. Therefore, the neighbor  of each vertex of $T_0$ is $v_{\frac{r-1}{2}}$ that is colored with color 2. It may happen (see {\bf [P8]}) that $F_{v_{\frac{r-1}{2}}}$ share an edge with face $F_{v_{\frac{r+1}{2}}}$, in such a case, the last two vertices in $T_0$ belong to these both faces and therefore they have both $v_{\frac{r-1}{2}}$  and $v_{\frac{r-1}{2}}$ as neighbors, but this is not a problem since $v_{\frac{r-1}{2}}$ is colored with color 1 (and vertices in $T_0$ are colored with 0.

In order to complete the coloring (the rest of vertices in $\mathbb{S}^2_{\text{\tiny N}}$ and all the vertices in $\mathbb{S}^2_{\text{\tiny S}}$), we need to partition $\mathbb{S}^2_{\text{\tiny N}}$ into 3 regions as follows:

\begin{itemize}
    \item $R_{0,1} := \text{ boarded by } \bar Q_0, \bar Q_1\text{ and the arc of the equator between } v_0 \text{ and } z_y,$
    \item $R_{1,2} := \text{ boarded by } \bar Q_1, \bar Q_2\text{ and the arc of the equator between } z_y \text{ and } z_x$
    and
    \item  $R_{0,2} := \text{ boarded by } \bar Q_0, \bar Q_2\text{ and the arc of the equator between } z_x \text{ and } v_0,$ see Figure \ref{DisjointPaths}
\end{itemize}

The goal of such a partition is to divide the set of faces lying in $\mathbb{S}^2_{\text{\tiny N}}$ into three parts (each of these parts is partitioned into faces). We may color the vertices of each region $R$ using only two colors and therefore the vertices of any face $f$ lying in $R$ would be colored either two colors (or three if $f$ contains the vertex $v$). Notice that $f$ is the dual face of a vertex $w$ lying in the equator or in $\mathbb{S}^2_{\text{\tiny S}}$. We would then always have a color left (other than color 3) to be used to color $w$. 

We shall color the vertices lying in the interior of region $R_{0,1}$ (similarly for the other two regions). Recall that all the vertices $v_i$ in the cycle $C_D$ are already well colored. Let $u$ be a vertex in the interior of $R_{0,1}$. Then,

$$c(u)=\left\{\begin{array}{ll}
c(v_i) & \text{ if } u\in P_{i} \text{ for some } i,\\
0 \text{ or } 1& \text{ otherwise}.
 \end{array}\right.$$

If $u\in P_{i}$ then $u$ and $v_i$ are both in the same dual face, say $F_{v_j}$, and therefore they are both neighbors of $v_j$ (which is a vertex in $C_D$). Since the $v_i$'s are well colored then $v_i$ is well colored with respect to $v_j$ and therefore any vertex of  $F_{v_j}$ (which is also a neighbor of $v_j$ in $\DDiag_G$) having the same color as $v_i$ would be well colored as well, in particular,  $u$ (which has the same color as $v_i$) is properly colored with respect to $v_j$.

If $u\not\in P_{i}$ then $u$ would belong to a face $F$ lying within the region $R_{0,1}$. Such vertices are colored with colors 0 or 1 (or 3 if the $F$ touches vertex $v$). Since $F$ is the dual face $F_{w}$ for some vertex $w$ lying in $\mathbb{S}^2_{\text{\tiny S}}$ then it would be enough to color $c(w)=2$. 


 On this way, we can always find a proper 4-coloring (with colors 0,1,2 and 3) in which $v$ is the only vertex having color 3, as desired. 
\end{Proof1}

Let us observe that if $R(V)$ is a Reuleaux polyhedron and $G$ its 1-skeleton,   then as we said before,  $G$ is an involutive polyhedral graph and $\DDiam_V\cong \DDiag_G$.  Therefore,  for every $v\in V$,  Lemma \ref{critical} induces a partition of $V$ in four parts,  such that each part has diameter less than $\ddiam (V)$ and one part contains only $v$.

\subsection{Critical partition of Reuleaux polyhedra}

 In this subsection, we will extend the``critical partition" of the vertices of any Reuleaux polyhedra $R(V)$ given by Lemma \ref{critical},  to the rest of the points in $R(V)$ (see Theorem \ref{MeissnerThm}).  To this end,  we will recall to  the reader some of the classical notions arising from Meissner polyhedrons. These definitions have been widely used since the first construction of Meissner and Schilling in 1912, (see \cite{martini2019bodies} and the references therein). Here, we will follow the careful description of the boundary of Reuleaux polyhedra given by Hynd \cite[Section 4.1,  pp88]{ hynd2024density} see also \cite{montejano2017meissner}.

Let $(e,e')$ be a pair of dual edges of a Reuleaux polyhedron $R(X)$ assume $e=\{a, b\}$ and $e'=\{c,  d\}$ then $|a-c|=|a-d|=|b-c|=|b-d|=1$, and $e'$ is a piece of the circle obtained by the intersection of the spheres with radius $1$ and centers at $a$ and $b$ (similarly $e$ is a piece of the circle obtained by the intersection of the spheres with radius $1$ and centers at $c$ and $d$).  

Let $H$ and $L$ be the planes generated by $a,c,d$ and $ b,c,d$ respectively.  Let $H^+$ and $L^+$ be the half-spaces that contains $e'$,  and let $W_e':=H^+\cap L^+\cap R(V)$  be the \emph{wedge} along the edge $e'$. It is  known, Lemma 4.2 of \cite{hynd2024density}) that a point $x\in W_e$ may only have distance greater or equal to $1$ with points in the wedge $W_{e'}$.  Recall, that since Reuleaux polyhedra are standard ball polyhedra, then their 1-skeleton is a polyhedral graph.  Therefore, if $x$ is not a vertex but a point lying on the boundary of $R(V)$ then it is on a facet $F_w$, that is dual to some vertex $w$ and it is surrounded by edges, say $e_i=\{a_i,b_i\}$, $i=1,\dots k$.  Furthermore if such $x$ does not lay in a wedge then it will be surrounded by the corresponding wedges $W_{e_i}$.



 

 \begin{theorem}\label{MeissnerThm}
     Let $R(V)$ be a Reuleaux polyhedron. For every $\epsilon>0$ and every $v\in V$, there is a partition $P_{\epsilon, v}=\{P_1, P_2, P_3, P_4\}$ of $R(V)$, such that $v\in P_1$, $\ddiam P_1\leq \epsilon$ and $\ddiam P_i <1$, for $i=2,3,4$.
 \end{theorem}

\begin{figure}
        \centering
    \begin{subfigure}[b]{0.3\textwidth}
    \centering
        \includegraphics[scale=0.5]{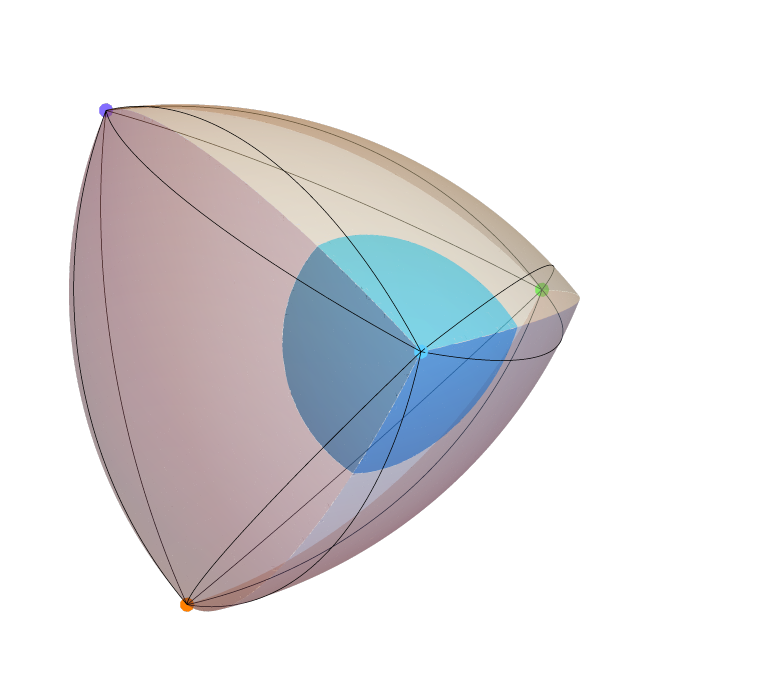}
        \caption{$P_1$ and vertices partition}
        \label{FigVerticesPartition}
    \end{subfigure}
     \hfill
        \centering
    \begin{subfigure}[b]{0.3\textwidth}
    \centering
        \includegraphics[scale=0.5]{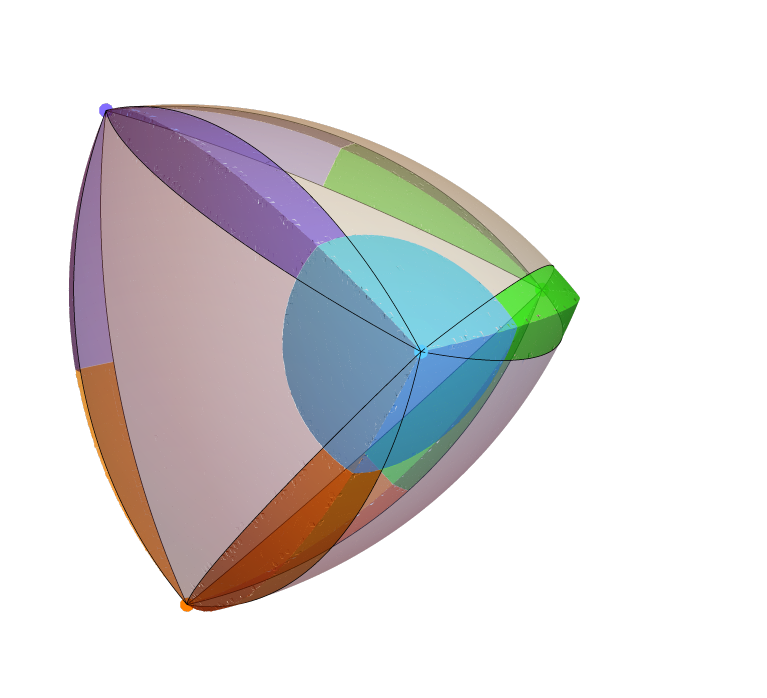}
        \caption{Wedges partition}
        \label{FigWedgesPartition}
    \end{subfigure}
     \hfill
        \centering
    \begin{subfigure}[b]{0.3\textwidth}
    \centering
        \includegraphics[scale=0.5]{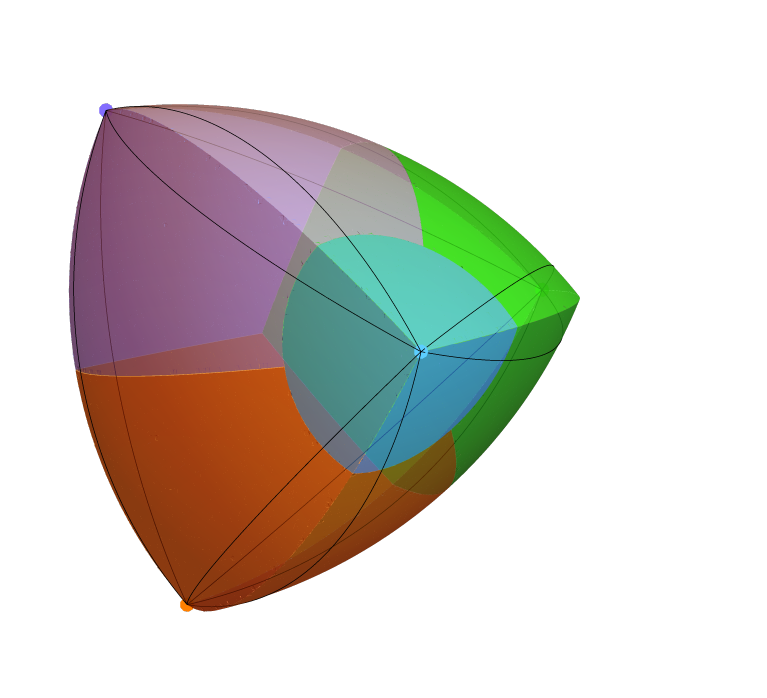}
        \caption{The whole partition}
        \label{FigWholePartition}
    \end{subfigure}
     \hfill
    \centering
    \caption{Critical partition of the tetrahedron}
    \label{TetraPartition}
   
\end{figure}

 \begin{proof}

Let $v$ be a vertex in $R(V)$ and $\epsilon >0$.  By Lemma \ref{critical}, we know that we can color the vertices $V$ of $R(V)$ with four colors, say $1, 2, 3$ and $4$, in such a way, that $v$ is the only vertex with color, say, $1$.  We will define a partition $P_{\epsilon, v}$ by extending this critical coloration to the rest of the Reuleaux polyhedron as follows.  We say that a vertex $a\in V$ is in part $P_a$ if its corresponding color is $i_a$ for some  $i_a\in \{1,2,3,4\}$, clearly if $a=v$, $P_a=P_v=1$ any other vertex $a\not=v$ will have color $i_a\in \{2,3,4\}$. 

Let $r$ be the minimum among the distances between $v$ and any vertex in $V\setminus \{v\}$.  Let $\epsilon_1:=\min\{\frac{r}{2}, \frac{\epsilon}{2}\}$ and let $P_1:=B_{\epsilon_1}(v)\cap R(V)$, where $B_{\epsilon_1}(v)$ is the open ball with radius $\epsilon_1$ and center at $v$. It is clear that $v$ is the only vertex in $P_1$ and $\ddiam P_1\leq \epsilon$. Denote by $\bar{P}_1$ the closure of 
$P_1$. 

    Now, we are going to partition $R(V)\setminus P_1$ (which is a compact set) in three parts, such that each part is a compact set with diameter less than $1$. This means that $\cup_{i=2,3,4}P_i=R(V)\setminus P_1$ and each $P_i$ is a compact set. Let us recall that this is not a partition in the classical way because $P_i\cap P_j$ might not be empty, but in order to prove that each part has diameter less than 1, it is easier to work with compact sets since all the diameters are going to be achieved in each part of the partition.
    
    Let us begin by extending the coloring of the set of vertices to every point of $(R(V)\setminus P_1)$ by following steps 1 to 4 in this  order.

     \begin{enumerate}
     

         \item \emph{Partitioning points that lay in a wedge.} \\ 
         Let $x\in R(V)\setminus P_1$ be a point in a wedge  $W_e$ generated by the  edge $e=\{a,  b\}$ with vertices $a, b$.\\ 
         We have two cases:
         \begin{itemize}
             \item[a)] If $a=v$ assign $x$ part $P_b$. (Similarly if $v=b$ assign $x$ part $P_a$ )
             \item[b)] If $a\not=v$ and  $b\not=v$ then 
             add $x$ to $P_a$ if $x$ is closer to $a$ then $b$,  otherwise $x$ will be in part $P_b$. (Clearly if the distances from $x$ to $a$ and to $b$ are the same,  then $x$ will be in both parts $P_a$ and $P_b$)
         \end{itemize}
         
\end{enumerate}

         Notice  that this partition induces compact sets in each wedge. Suppose $e'=\{c,  d\}$ is the dual edge of $e=\{a,  b\}$. Since $c$ and $d$ cannot be in the same partition containing neither $a$ nor $b$,  then,  this partition of the wedges is valid by compactness.  We may thus say that every point on a wedge  $x\in W_e$ has been added to partition $P_x$ (where $P_x=P_a$ of $P_x=P_b$) in such a way that all the points belonging to the same partition have diameter less than $1$ (so far so good,  see Figure \ref{FigWedgesPartition}).

     Let $x$ be on the boundary of  of $R(V)\setminus P_1$.   Observe that if $x$ is not a vertex and is not on a wedge then two possibilities may occur: 
        
        \begin{itemize}
            \item $x$ is in a face of $R(V)$ but not in a wedge.       
            \item $x\in \bar{P_1}\cap (R(V)\setminus P_1)$.
        \end{itemize} 
        
\begin{enumerate}[resume]
 \item \emph{Partitioning points that are in the boundary of $R(V)$ but not in the interior of a wedge.} \\
If  $x$ is in the boundary of $R(V)$, then $x\in F_w\setminus P_1$ where $F_w$ is the dual facet of some vertex $w\in V$. As we observed above $x$ is surrounded by some wedges $W_{e_i}$, for $i=1,\dots ,k$. Let $y$ be a point in one of these wedges say, $W_{e_j}$ such that $y$ is the closest point to $x$.  We add $x$ to partition $P_y$. Since $F_w$ is dual to $w$ then $|w-a_j|=1=|w-b_j|$ and therefore by Lemma \ref{critical} $P_w\not= P_{a_j}$ and $P_w\not= P_{b_j}$, hence $P_y\not= P_w$. Furthermore, $w$ is the only point at distance one from $x$, and so any point in $P_y$ is at distance less than one from $x$, as needed.

\item \emph{ Partitioning points that are in $ \bar{P_1}\cap (R(V)\setminus P_1)$. }
Let $x\in \bar{P_1}\cap (R(V)\setminus P_1)$.  Since $x$ is not on a wedge then $x$ and is not at the boundary of $R(V)$ and so  $x$ is in the interior of $R(V)$.  We thus add $x$ to partition $P_y$,  where $y$ is the closest point to $x$ on the boundary of $R(V)$.




         \item \emph{Partitioning points that are in the interior of $R(V)\setminus P_1$.} 
         Let $c$ be the circumcenter of $R(V)\setminus P_1$. 
         Let $x\not=c$ be an interior point of $R(V)\setminus P_1$ and let $l_{xc}$ be the half-line starting at $c$ and passing through $x$. We add the point $x$ to partition $P_w$ where $w$ is the point of intersection of $l_{xc}$ and the boundary of $R(V)\setminus P_1$.

         It is known that the maximum diameter of a Reuleaux polyhedron is $\sqrt{3}-\sqrt{2}/2\approx 1.02$ (see \cite{martini2019bodies}).  Then the circumsphere of any Reuleaux polyhedron have at most $(\sqrt{3}-\sqrt{2}/2) (\sqrt{\frac{3}{8}})<1$ due to the classic Jung's Theorem (see \cite{martini2019bodies}[Theorem 15.2.1] or the original version \cite{jung1899kleinste}).    Therefore the circumradius of $R(V)\setminus P_1$ must be less than $1$, which implies that the interior of $R(V)\setminus P_1$ has been well partitioned.

     \end{enumerate}
 \end{proof}
 
 In Figure \ref{TetraPartition},   we represent a possible critical partition for the tetrahedron  step by step,  where the blue color represents $P_1$.
 
\section{Main Results}\label{Main Results}
In this section, we prove our main contributions. We first show the validity of Conjecture \ref{vazcritical} (see Theorem \ref{Conjballteo}), which leads us to the proof of our main result. Then towards the end of this section, we present a special configuration of points that is critical but not strongly critical for the Vázsonyi problem.

\subsection{Reuleaux polyhedra in the Vázsonyi problem}\label{ReuleauxSection}
In order to show Conjecture \ref{vazcritical}, we need the following

\begin{lemma}\label{3connected}
    Let $V\subset\mathbb{R}^3$ be an extremal configuration for the Vázsonyi problem. Then, the 1-skeleton of 
    $\mathcal{B}(V)$ is planar,  simple and $3$-connected if and only if $V$ is strongly critical. 
\end{lemma}

\begin{proof} 

Let us denote by $G$ the 1-skeleton of $\mathcal{SF(B}(V))$.

(\emph{Necessity}) Suppose that $V$ is strongly critical, that is, $V$ does not have an extremal proper subset. Since in particular $V$ is an extremal configuration, by Theorem \ref{dualedges}, $G$ admits a canonical involution, say $\varphi$. Furthermore, by Theorem \ref{GHS} $V$ is tight and by Theorem \ref{2connected}, $G$ is a 2-connected planar graph.

It is known \cite{bezdek2007ball}[pp-19] that $G$ has no loops and, it was shown \cite{hynd2024density}[Theorem 3.6] that $G$ has no multiple edges when $V$ is extremal for the Vázsonyi problem.  Therefore, $G$ is simple.

We shall prove now that $G$  is 3-connected. Let us proceed by contradiction.  Suppose that $G$ admits a 2-cutting set,  say $\{x,  y\}$.  Let $A_1,\dots ,A_k$, $k\geq 2$, be the connected components of $G\setminus \{x,y\}$. Clearly,  $x$ and $y$ are both adjacent to each $A_i$. 

Let $B_i:=\varphi(V(A_i))$,  this is $$B_i=\bigcup_{a\in V(A_i)} \varphi (a),$$  for each $1\leq i \leq k$.  Let us observe that $B_i$ is a subgraph of $G$ and the union of cycles (the dual faces of $V(A_i)$). 

We observe the following claims: 
 \smallskip

\begin{itemize}

   \item[\textbf{[a]}] \label{NotSharing} \emph{$B_i$ and $B_j$  share no edges,  for all $i\neq j$.  Furthermore $F_x$ and $F_y$  share at least one edge with each $B_i$}
   
   Let $e=F_v\cap F_w$ be an edge of $B_i$,  with $v,w\in V(G)$.  Then either $v$ or $w$ have to be in $V(A_i)$. Without loss of generality,  suppose $v\in V(A_i)$. 
   
   If we suppose that $e$ is in some $B_j$,  with $i\neq j$,  then $w$ must be in $A_j$.  Since the dual edge of $e$ is $e'=\{v, w\}$,  then $A_i$ and $A_j$  would be connected  by this edge,  which is a contradiction. 
   
   Then $e$ is not in $B_j$,  which means that $w$ is either in $A_i$ or $w\in \{x,  y\}$.  In this case,  we say that $e$ is a \emph{boundary } edge.  Furthermore,  since $\{x, y\}$ is a cutting set,  then for $x$ and $y$,  there is at least vertices $v_1,v_2\in V(A_i)$ (not necessarily different) such that $x$ is adjacent to $v_1$ and $y$ is adjacent to $v_2$.  Therefore,  $e_1=\{x,  v_1\}$ and $e_2=\{y,  v_2\}$ are boundary edges and their dual edges $e_1'$ and $e_2'$ share an edge with $F_x$ and $F_y$ respectively.
   
        \item[\textbf{[b]}] \emph{$B_i$  has more than three vertices.  }
   
   Since  $A_i$ has more than two vertices  (otherwise,  $A_i$ would consist of a dangling vertex which is not possible since  $V$ is strongly critical). Let $u,v\in V(A_i)$,  then $F_u$ and $F_v$ have at least three vertices each (because $V$ is strongly critical) and  since $F_u$ and $F_v$ are different cycles,  then $V(F_u\cup  F_v)\geq 4$.   
   
    \item[\textbf{[c]}]\label{BoundaryConnected} \emph{$B_i$ is a connected graph. }
    
      Let us prove that ${B}_i$ is connected for each $1\le i\le k$. Indeed, Let $p,q\in V({B}_i)$ we show that there is a path $\gamma_{p,q}$ joining $p$ and $q$. Suppose that $p\in F_{r}$ and $q\in F_{s}$ were $F_{r}$ and $F_{s}$ are some faces in ${B}_i$, $r,s\in A_i$. Since $A_i$ is connected then there exists a path $\gamma[r,s]$ between the vertices $r$ and $s$. Assume first that $\gamma[r,s]$ consists of one edge. Then, $F_r$ and $F_s$ must share one edge.  We can thus construct a path from $p$ to $q$ by a proper sequence of vertices in $F_r$ and $F_s$. Now, we can clearly proceed by induction on the length of $\gamma[r,s]$ when it is greater than or equal to 2.
      
          \item[\textbf{[d]}] \emph{There is an embedding of $G$ in the plane such that $F_x$ and $F_y$  are contained in the outer face of $B_i$.  Furthermore the boundary of this outer face is a cycle.}
      
Let $\hat{G}$ be the embedding of $G$ in the plane such that $F_x$ is the outer face of $\hat{G}$.
     
     If $F_y$ is not in the outer face of $B_i$ for this embedding,  then $F_y$  must be totally contained in a face $F$ of $B_i$ different than its outer face (see Figure \ref{FigLemma2Nueva2}).
     
     Since $F$ is a cycle and all its edges are also edges of $F_y$ by \textbf{[a]},  then $F$ must be exactly $F_y$. Which is a contradiction with \textbf{[a]} because any other $B_j$ would not be able to share an edge with $F_y$.
     
     Therefore,  $F_x$ and $F_y$  are in the outer face of $B_i$ in $\hat{G}$. Let us observe that this fact implies that all the dual faces of vertices not in $V(A_i)$ also have to be in this outer face  of $B_i$,  in other words,  all the faces of $B_i$ except its outer face in $\hat{G}$ are dual faces of the vertices of $A_i$.
     
     Finally,  let us recall that a graph is \emph{nonseparable} if there is not a decomposition of the graph  into two nonempty connected subgraphs which have just one vertex in common.  Then $B_i$ is a nonseparble graph because is the union of cycles (the dual faces of $V(A_i)$) by edges.   In \cite{bondy2008graph}[Theorem 10.7],  the authors prove that all the faces of a planar nonseparable graph other than the complete graphs with one and two vertices $K_1$ and $K_2$ are cycles.  But $B_i$ is neither $K_1$ nor $K_2$ because it has more than three vertices by \textbf{[b]},  the outer face of $B_i$ is a cycle.
     
     Let us denote the edges of this outer face by $\partial B_i$ .
     
      \begin{figure}[ht]
    \centering
    \begin{subfigure}[b]{0.55\textwidth}
        \includegraphics[width=\textwidth]{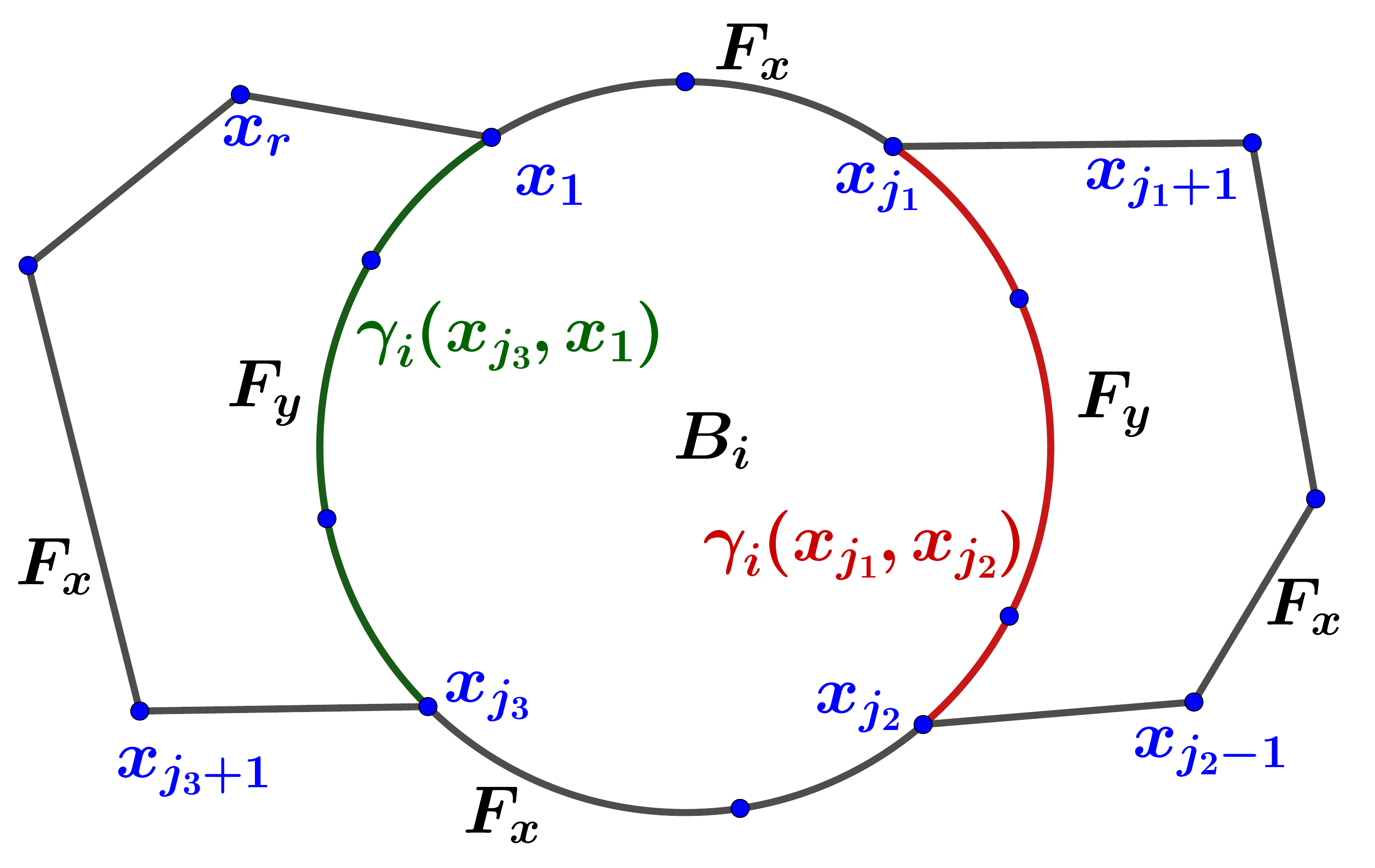}
        \caption{}
        \label{FigLemma2Nueva1}
    \end{subfigure}
     \hfill
    \begin{subfigure}[b]{0.35\textwidth}
        \includegraphics[width=\textwidth]{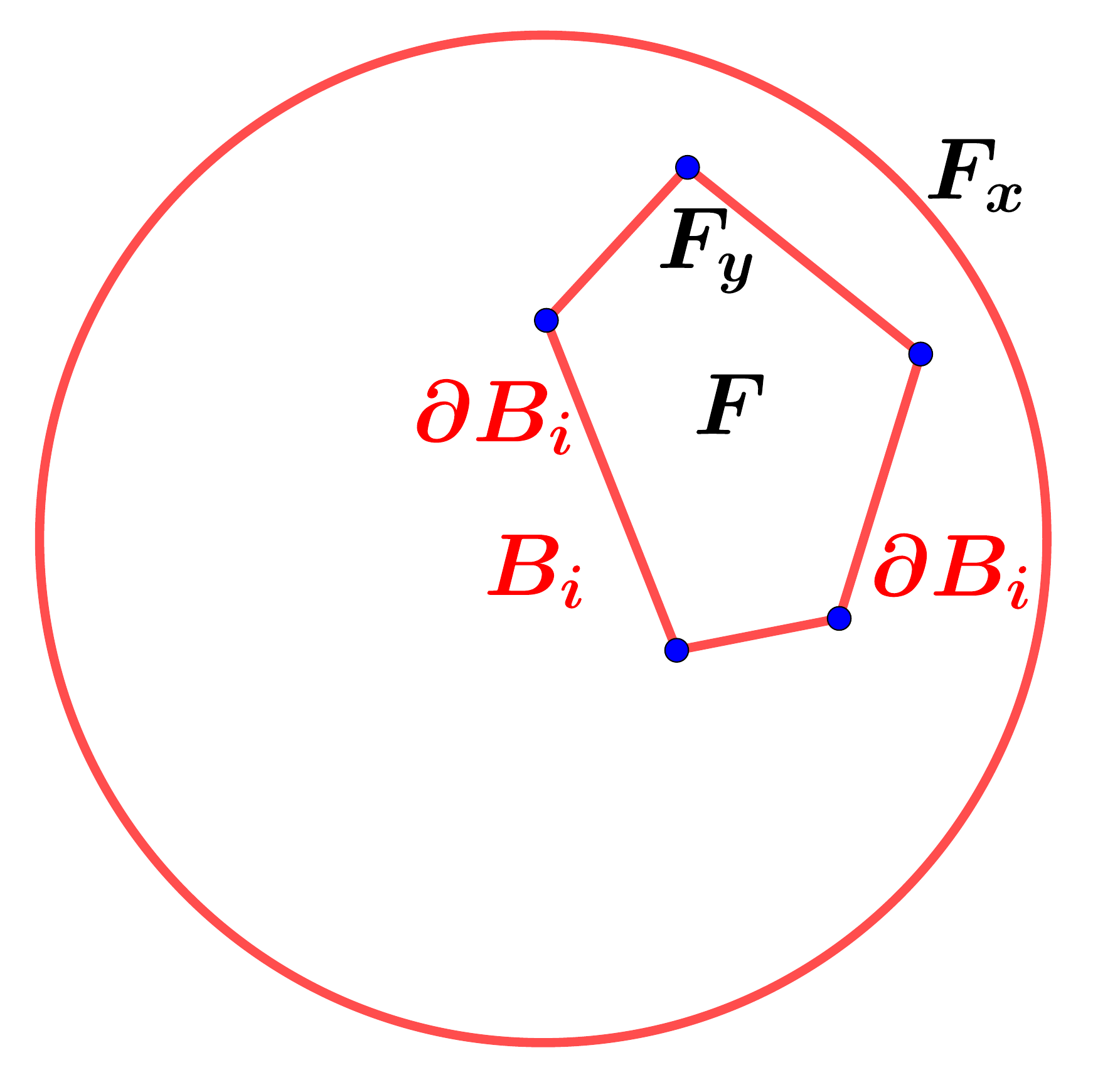}
        \caption{}
        \label{FigLemma2Nueva2}
    \end{subfigure}
    
    \caption{}
    \label{FigLemma2Nueva}
   
\end{figure}
   
   \item[\textbf{[e]}]\label{FxFyPaths} \emph{$F_x$ (respectively $F_y$) is the union of the $k$ paths $\gamma_{i}^x$, (respectively $\gamma_{i}^y$), and possibly an extra edge shared by $F_x$ and $F_y$}
   
  Let us first prove the statement for $F_x$.  An analogous procedure can be done for $F_y$.  Since $F_x$ is a cycle,  let $\{x_1,...,  x_r\}$ be its set of vertices.
  
  Let $N(x)\subset V(G)$ be the neighbours of $x$,  that is,  the set of all the vertices adjacent to  $x$.  By the involution $\varphi $,   we know that all the edges of $F_x$  are of the form $\{x,  v\}'$  for some $v\in N(x)$.  Furthermore,  for all $v\in N(x)$  not equal to $y$,   there is a $A_i$ such that $v\in A_i$,  thus $\{x,  v\}'$ is also an edge of $B_i$.  Therefore,  if $y\in V_x$,  all the edges of $F_x$ are in some $B_i$  except the edge $\{x, y\}$,  which is the unique edge shared by $F_x$ and $F_y$.  If not,  all the edges of $F_x$ are in some $B_i$.  
   
 	By \textbf{[a]},  we know that $F_x$ and $B_i$ share at least one edge  for all $i=1,...,k$.  So it remains to prove that they share exactly one path.  Let us proceed by contradiction,  without loss of generality,  suppose that  $x_1,...x_{j_1}$ and $x_{j_2},..., x_{j_3}$  are two disjoint paths that share edges with $B_i$,  with $1< j_1<j_2<j_3<r$,  but the paths $x_{j_1},...,  x_{j_2}$ and $x_{j_3},..., x_r,x_1$ do not share edges with $B_i$. 
 	
 	Let us observe that $x_{j_1}$ and $x_{j_2}$  are vertices in $\partial B_i$, then,  by \textbf{[d]}, there is a path $\gamma_i (x_{j_1},  x_{j_2})$ totally contained  in $\partial B_i$ connecting $x_{j_1}$ and $x_{j_2}$ such that the cycle $C$ limited by $x_{j_1},...,  x_{j_2}$ and $\gamma_i (x_{j_1},  x_{j_2})$ does not contain $B_i$ (see Figure \ref{FigLemma2Nueva1}).  Since  $\gamma_i (x_{j_1},  x_{j_2})$ must share all its edges with $F_y$ by \textbf{[a]},  then $F_y$  must be contained in $C$.  Similarly,  $F_y$  must be contained in the cycle limited by $x_{j_3},..., x_r,x_1$ and  $\gamma_i (x_{j_3},  x_{1})$.  Since these cycles no share interior,  then $F_y$ cannot be in both cycles at the same time,  holding a contradiction. 
 	
 	Therefore,  there is a unique path $\gamma_i^x $ shared by $F_x$  and $B_i$,  for all $i=1,...,k$.  As we said at the beginning,  we may repeat this process for $F_y$,   then there is a unique path $\gamma_i^y $ shared by $F_y$  and $B_i$,  for all $i=1,...,k$.
 	
 	\item [\textbf{[f]}] \emph{$V(F_x)\cap V(F_y)\neq \emptyset$ (see Figure \ref{FxFy}). }
 	
 	Let us recall that by \textbf{[e]}  all the edges of $\partial B_i$ are in the paths $\gamma_i^x $ and $\gamma_i^y $ and by \textbf{[d]}, $\partial B_i$ is a cycle.   Then  $\gamma_i^x $ and $\gamma_i^y $ must share final vertices in order to create $\partial B_i$.  
 	
 	Let $w_i$ and $z_i$ be those end points of $\gamma_i^y\subset F_y$ and $\gamma_{i}^x\subset F_x$. Then $w_i,  z_i\in V(F_x)\cap V(F_y)$.

   \end{itemize}

\begin{figure}[ht]
    \centering
    \begin{subfigure}[b]{0.4\textwidth}
        \centering
        \definecolor{qqqqff}{rgb}{0.,0.,1.}
        \begin{tikzpicture}[line cap=round,line join=round,>=triangle 45,x=1.0cm,y=1.0cm, scale=2]
        \clip(1.25402718810581,0.213276119392723) rectangle (4.350271017830777,2.638947957841113);
       \draw [shift={(2.858159908529994,1.5397537865988118)},line width=2.pt]  plot[domain=0.4510082617889014:2.649314860432547,variable=\t]({1.*0.9737889943711315*cos(\t r)+0.*0.9737889943711315*sin(\t r)},{0.*0.9737889943711315*cos(\t r)+1.*0.9737889943711315*sin(\t r)});
        \draw [shift={(2.8779916860254047,2.500703039491826)},line width=2.pt]  plot[domain=3.6598741613041152:5.723634268096919,variable=\t]({1.*1.0107289124617371*cos(\t r)+0.*1.0107289124617371*sin(\t r)},{0.*1.0107289124617371*cos(\t r)+1.*1.0107289124617371*sin(\t r)});
        \draw [shift={(5.618765535888748,0.9493968436260437)},line width=2.pt]  plot[domain=2.6475513588296575:3.3990389044636173,variable=\t]({1.*2.1400918835376124*cos(\t r)+0.*2.1400918835376124*sin(\t r)},{0.*2.1400918835376124*cos(\t r)+1.*2.1400918835376124*sin(\t r)});
        \draw [shift={(3.4448782070973256,1.2077686982774836)},line width=2.pt]  plot[domain=-1.4416421202170566:1.2050470763307461,variable=\t]({1.*0.8100109315716193*cos(\t r)+0.*0.8100109315716193*sin(\t r)},{0.*0.8100109315716193*cos(\t r)+1.*0.8100109315716193*sin(\t r)});
        \draw [shift={(2.3377927127863902,1.3504477178919116)},line width=2.pt]  plot[domain=2.8168117080027404:5.027461264368392,variable=\t]({1.*1.0419815591154444*cos(\t r)+0.*1.0419815591154444*sin(\t r)},{0.*1.0419815591154444*cos(\t r)+1.*1.0419815591154444*sin(\t r)});
        \draw [shift={(1.30436845609526,0.32700823662774886)},line width=2.pt]  plot[domain=0.02414193454335175:1.536945730648194,variable=\t]({1.*1.3567142679317086*cos(\t r)+0.*1.3567142679317086*sin(\t r)},{0.*1.3567142679317086*cos(\t r)+1.*1.3567142679317086*sin(\t r)});
        \draw [line width=2.pt] (1.3502852726306116,1.682945274535015)-- (2.,2.);
        \draw (2.7844613791934902,2.2336054658063067) node[anchor=north west] {$B_1$};
        \draw (3.73418908235957,1.3669789366672604) node[anchor=north west] {$B_2$};
        \draw (1.7219535112764384,1.2482629737715005) node[anchor=north west] {$B_k$};
        \draw [color=qqqqff](2.802268773627854,1.2541987719162886) node[anchor=north west] {$F_x$};
        \draw [color=qqqqff](1.4726499891953424,2.245477062095883) node[anchor=north west] {$F_y$};
        \draw (2.600451636705062,1.472860275445264) node[anchor=north west] {$\gamma_{1}^x$};
        \draw (2.10184459254287,2.595689152638374) node[anchor=north west] {$\gamma_{1}^y$};
        \draw (2.0840371981085064,2.162375888068851) node[anchor=north west] {$w_1$};
        \draw (3.4967571565680498,2.126761099200123) node[anchor=north west] {$z_1$};
        \begin{scriptsize}
        \draw [fill=qqqqff] (2,2) circle (1.5pt);
        \draw [fill=qqqqff] (3.734577877138053,1.9642023109916569) circle (1.5pt);
        \draw [fill=qqqqff] (3.02,1.5) circle (1.5pt);
        \draw [fill=qqqqff] (3.47888966217747,0.9798026833934106) circle (1.5pt);
        \draw [fill=qqqqff] (3.54920392129163,0.4045041997320978) circle (1.5pt);
        \draw [fill=qqqqff] (4.207601074815133,0.9350572457753085) circle (1.5pt);
        \draw [fill=qqqqff] (2.660687374303602,0.3597587621139957) circle (1.5pt);
        \draw [fill=qqqqff] (1.4525605586148451,0.8008209329210022) circle (1.5pt);
        \draw [fill=qqqqff] (1.3502852726306116,1.682945274535015) circle (1.5pt);
        \draw [fill=qqqqff] (2.0981733013903185,1.4272570595744316) circle (1.5pt);
        \draw [fill=black] (3.2487702687129443,0.3661509674880103) circle (0.5pt);
        \draw [fill=black] (3.1337105719806817,0.3597587621139957) circle (0.5pt);
        \draw [fill=black] (3.018650875248419,0.3597587621139957) circle (0.5pt);
        \end{scriptsize}
        \end{tikzpicture}
        \caption{The boundaries of $F_x$, $F_y$ and ${B_i's}$}
        \label{FxFy}
    \end{subfigure}
    \hfill
    \begin{subfigure}[b]{0.4\textwidth}
        \centering
        \definecolor{ududff}{rgb}{0.30196078431372547,0.30196078431372547,1.}
        \definecolor{qqwuqq}{rgb}{0.,0.39215686274509803,0.}
        \definecolor{ffqqqq}{rgb}{1.,0.,0.}
        \definecolor{qqqqff}{rgb}{0.,0.,1.}
        \begin{tikzpicture}[line cap=round,line join=round,>=triangle 45,x=1.0cm,y=1.0cm]
        \clip(-0.3545183746673877,-0.5877784135725105) rectangle (6.4,4.726977648084146);
        \fill[line width=2.pt,color=qqwuqq,fill=qqwuqq,fill opacity=0.03999999910593033] (-0.1499574516683648,4.143068538983268) -- (-0.1625796070626519,-0.1295310619829169) -- (6.287341799418058,-0.14215321737720402) -- (6.262097488629483,4.187246082863274) -- cycle;
        \draw [shift={(2.5,2.)},line width=2.pt]  plot[domain=2.214297435588181:4.068887871591405,variable=\t]({1.*2.5*cos(\t r)+0.*2.5*sin(\t r)},{0.*2.5*cos(\t r)+1.*2.5*sin(\t r)});
        \draw [shift={(-0.5,2.)},line width=2.pt]  plot[domain=-0.9272952180016123:0.9272952180016122,variable=\t]({1.*2.5*cos(\t r)+0.*2.5*sin(\t r)},{0.*2.5*cos(\t r)+1.*2.5*sin(\t r)});
        \draw [shift={(6.5,2.)},line width=2.pt]  plot[domain=2.214297435588181:4.068887871591405,variable=\t]({1.*2.5*cos(\t r)+0.*2.5*sin(\t r)},{0.*2.5*cos(\t r)+1.*2.5*sin(\t r)});
        \draw [shift={(3.5,2.)},line width=2.pt]  plot[domain=-0.9272952180016123:0.9272952180016122,variable=\t]({1.*2.5*cos(\t r)+0.*2.5*sin(\t r)},{0.*2.5*cos(\t r)+1.*2.5*sin(\t r)});
        \draw [line width=2.pt,color=ffqqqq] (1.,4.)-- (5.,4.);
        \draw [line width=2.pt,color=ffqqqq] (1.,4.)-- (5.,0.);
        \draw [line width=2.pt,color=ffqqqq] (1.,0.)-- (5.,4.);
        \draw [line width=2.pt,color=ffqqqq] (1.,0.)-- (5.,0.);
        \draw [line width=2.pt,color=ffqqqq] (1.200613175520355,2.8997862326459884)-- (5.271258290177946,2.8240533002802657);
        \draw [line width=2.pt,color=ffqqqq] (0.42435061877169833,2.218189841354485)-- (5.258636134783659,2.2497452298402028);
        \draw [line width=2.pt,color=ffqqqq] (1.2069242532174986,0.9749075350172048)-- (5.176592124720793,1.012774001200066);
        \draw (0.474584849463441,4.200246902919388) node[anchor=north west] {$x$};
        \draw (0.4769484454498973,0.28695321989020882) node[anchor=north west] {$y$};
        \draw (5.1571540054135925,4.169701286865213) node[anchor=north west] {$w$};
        \draw (5.141881197386505,0.24458962390375252) node[anchor=north west] {$z$};
        \draw (0.9189497778968394,2.0696901831406964) node[anchor=north west] {$A_1$};
        \draw (4.897516268953106,2.031508163072978) node[anchor=north west] {$B^-_1$};
        \draw [line width=2.pt,color=qqwuqq] (-0.1499574516683648,4.143068538983268)-- (-0.1625796070626519,-0.1295310619829169);
        \draw [line width=2.pt,color=qqwuqq] (-0.1625796070626519,-0.1295310619829169)-- (6.287341799418058,-0.14215321737720402);
        \draw [line width=2.pt,color=qqwuqq] (6.287341799418058,-0.14215321737720402)-- (6.262097488629483,4.187246082863274);
        \draw [line width=2.pt,color=qqwuqq] (6.262097488629483,4.187246082863274)-- (-0.1499574516683648,4.143068538983268);
        \draw [line width=2.pt] (1.,4.)-- (1.4229524427907232,4.5820671035965725);
        \draw [line width=2.pt] (1.,4.)-- (0.5829480013009162,4.551521487542398);
        \draw [line width=2.pt] (5.,4.)-- (5.004425925142718,4.475157447406961);
        \draw [line width=2.pt] (1.,0.)-- (0.5524023852467416,-0.3968683132339181);
        \draw [line width=2.pt] (1.,0.)-- (1.0334958380999948,-0.43505033330163656);
        \draw [line width=2.pt] (1.,0.)-- (1.6902265832647527,-0.38159550520683067);
        \draw [line width=2.pt] (5.,0.)-- (4.637878532492621,-0.4045047172474618);
        \draw [line width=2.pt] (5.,0.)-- (5.6000654381991275,-0.3892319092203744);
        \draw [color=qqwuqq](-0.105490508312226507,3.9973354701338535) node[anchor=north west] {$G_1$};
        \draw (2.9884152655671823,4.6820671035965725) node[anchor=north west] {$G_2$};
        \draw (3.0342336896484445,-0.10668496071925757) node[anchor=north west] {$G_2$};
        \begin{scriptsize}
        \draw [fill=qqqqff] (1.,0.) circle (2.5pt);
        \draw [fill=qqqqff] (1.,4.) circle (2.5pt);
        \draw [fill=qqqqff] (5.,0.) circle (2.5pt);
        \draw [fill=qqqqff] (5.,4.) circle (2.5pt);
        \draw [fill=qqqqff] (1.200613175520355,2.8997862326459884) circle (2.5pt);
        \draw [fill=qqqqff] (5.271258290177946,2.8240533002802657) circle (2.5pt);
        \draw [fill=qqqqff] (0.42435061877169833,2.218189841354485) circle (2.5pt);
        \draw [fill=qqqqff] (5.258636134783659,2.2497452298402028) circle (2.5pt);
        \draw [fill=qqqqff] (1.2069242532174986,0.9749075350172048) circle (2.5pt);
        \draw [fill=qqqqff] (5.176592124720793,1.012774001200066) circle (2.5pt);
        \draw [fill=ududff] (1.4229524427907232,4.5820671035965725) circle (2.5pt);
        \draw [fill=ududff] (0.5829480013009162,4.551521487542398) circle (2.5pt);
        \draw [fill=ududff] (5.004425925142718,4.475157447406961) circle (2.5pt);
        \draw [fill=ududff] (0.5524023852467416,-0.3968683132339181) circle (2.5pt);
        \draw [fill=ududff] (1.0334958380999948,-0.43505033330163656) circle (2.5pt);
        \draw [fill=ududff] (1.6902265832647527,-0.38159550520683067) circle (2.5pt);
        \draw [fill=ududff] (4.637878532492621,-0.4045047172474618) circle (2.5pt);
        \draw [fill=ududff] (5.6000654381991275,-0.3892319092203744) circle (2.5pt);
        \end{scriptsize}
        \end{tikzpicture}
        \caption{$A_1$, $B^-_1$, $G_1$ and $G_2$.}
        \label{ABG}
    \end{subfigure}
 
    \caption{}
   \label{VazFig}
\end{figure}

\begin{itemize}
   
    \item[\textbf{[g]}] \emph{For every $i\in \{1,\dots k\}$,  we have that $\{w_i,z_i\}$ is also a $2$-cutting set of $G$}\\
    We focus our attention to one of the connected component, say, $A_1$ and its image $B_1$. Let $A_1^+=A_1\cup\{x,y\}$ and $B_1^-=B_1\setminus \{w_1,z_1\}$. 
    Observe that any vertex in ${B}_1^-$ cannot be connected by a path to any other vertex in $B_i$, $i\neq 1$, in other words, $\{w_1,z_1\}$ is also a $2$-cutting set of $G$.  Furthermore observe that $B_1^-$ has the same ``shape" as $A_1$, that is, $B_1^-$ is connected, and the vertices $w_1$ and $z_1$ play the same role as the vertices $x$ and $y$ for $A_1$. Since $\varphi$ is involutive and $\varphi (A_1)={B}_1$ then $\varphi (B_1^-)={A_1^+}$. Thus we have that $A_1^+$ has the same shape as $B_1$ (see Figure \ref{ABG}). 

    \item[\textbf{[h]}] \emph{$\{x,y,w_1,z_1\}$ are four different vertices, moreover $w_1,z_1 \notin V(A_1)$} \\
    By the involutive properties of $\varphi$ we know that $x\notin F_x$ and $y\notin F_y$ and since $w_1$ and $z_1$ are in $F_x\cap F_y$ then $\{x,y\}\not=\{w_1,z_1\}$.  Now suppose that $w_1\in V(A_1)$ (the other case is similar), Then, $F_{w_1}\subset {B}_1$, now since $\{w_1,z_1\}\in F_{x}\cap F_{y}$, then $x,y\in F_{w_1}\cap F_{z_1}$, therefore $F_{z_1}\subset B_1$, which implies that $A_1\subset B_1$. Then $B_1=\varphi(A_1)\subset \varphi(B_1)= A_1$, implying that $A_1=B_1$ and therefore $\{x,y\}=\{w_1,z_1\}$ which leads to a contradiction.  


    \item[\textbf{[i]}] ${V(A_1^+)}\cap V({B_1})=\emptyset$\\
    We proceed by contradiction, suppose that there is $v\in V({A_1^+})\cap V({B_1})$, then there is path $\gamma[v,x]$ (completely contained in ${A}_1$)  joining $v$ to $x$. Since $w,z\notin V(A_1)\subset V(A_1^+)$ then $\gamma[v,x]$  contains neither $w_1$ nor $z_1$.  Any path starting from a vertex in $B_1^-$ that does not use either $z_1$ or $w_1$ must contain only vertices in $B_1$. The latter implies that $x\in B_1$,  which leads to a contradiction. .

\end{itemize}

We now count the number of diameters induced by $V(G)$.  Let $G_1$ be the subgraph generated by $V(A_1^+ \cup B_1)$ and let $G_2$ be the subgraph generated by  $V(G)\setminus V(A_1 \cup B_1^-)$ (see Figure \ref{ABG}). Since ${V(A_1^+)}\cap V({B_1})=\emptyset$ then $\{x,y,w_1,z_1\}=V(G_1 \cap G_2)$, so $$|V(G)|= |V(G_1)|+|V(G_2)|-4.$$

We have that $|E(\DDiam_{G})|=|E(\DDiam_{G_1})|+|E(\DDiam_{G_2})|-r $ where $r$ denotes the number of diameters having ends in $\{x,y,w,z\}$. Clearly, since $w_1,z_1\in F_x$ and $w_1,z_1\in F_y$,  the pairs $\{x,w_1\}, \{x,z_1\},\{y,w_1\}$ and  $\{y,z_1\}$ are diameters, thus  $r\geq 4$. 

Since $V$ is an extremal configuration, then 

$$|E(\DDiam_{G_1})|+|E(\DDiam_{G_2})|-r =|E(\DDiam_{G})|  =2|V(G)|-2=2(|V(G_1)|+2|V(G_2)|-4)-2,$$

and thus, 
\begin{equation} \label{eq1}
|E(\DDiam_{G_1})|+|E(\DDiam_{G_2})|=2(|V(G_1)|+2|V(G_2)|-10+r.
\end{equation}

Since $V$ is strongly critical then 
\begin{equation} \label{eq2}
|E(\DDiam_{G_1})|\leq 2|V(G_1)|-3 \text{ and } |E(\DDiam_{G_2})|\leq 2|V(G_2)|-3, 
\end{equation}

and thus, by adding these inequalities, we obtain

\begin{equation} \label{eq3}
  |E(\DDiam_{G_1})|+|E(\DDiam_{G_2})|\leq 2|V(G_1)|+2|V(G_2)|-6.
\end{equation}

By combining \eqref{eq1} with \eqref{eq3}, we have that $r=4$, that is,

\begin{equation} \label{eq4}
|E(\DDiam_{G_1})|+|E(\DDiam_{G_2})|=2(|V(G_1)|+2|V(G_2)|-6
\end{equation}

and so $\{x,w_1\}, \{x,z_1\},\{y,w_1\},\{y,z_1\}$ are the only diameters of $\DDiam_G$ on the set $\{x,y,z_1,w_1\}$. Furthermore, by combining \eqref{eq2} with \eqref{eq4}, we obtain that $|E(\DDiam_{G_1})|= 2|V(G_1)|-3$ and $|E(\DDiam_{G_2})|= 2|V(G_2)|-3$ are both odd integers.

We claim that $|E(\DDiam_{G_1})|$ is also an even integer, leading to the desired contradiction.
To this end, we first count the edges in $E(G_1)$ not having both ends in $\{x,y,w,z\}$, we denote by $\tilde{E}(G_1)$ such a set of edges. We know that, by construction,  the dual edge of an edge adjacent to a vertex $a\in V(A_1)$ is an edge in ${B_1}$ and, symmetrically,  the dual edge of an edge adjacent to a vertex $b\in V(B_1^-)$ is an edge in ${A_1^+}$
In other words, any edge in $\tilde{E}(G_1)$ will have its duals in $\tilde{E}(G_1)$. Then, the number of edges in $\tilde{E}(G_1)$ is even.

Now, we clearly have that

\begin{equation}\label{thecounting}
 \sum_{v\in A}\delta(v)+ \sum_{v\in B}\delta(v) +\sum_{v\in \{x,y\}}\restr{\delta(v)}{A}  +  \sum_{v\in \{w,z\}}\restr{ \delta(v)}{B}=2|\tilde{E}( G_1)|   
\end{equation}

where $\delta(v)$ denotes the degree of a vertex $v$ in $G$ and $\restr{\delta(v)}{S}$ the degree of vertex $v$ with endpoints only on set $S$. 

We observe that, by duality, the degree of each vertex $v\in A_1$ is the same as the number of vertices of its dual face and thus the number of diameters adjacent to $v$. Then, the diameters with one end in $A_1$ is $\sum_{v\in A_1}\delta(v)$. By the same argument, $\sum_{v\in B^-_1}\delta(v)$ gives the diameters with one end in $B^-_1$.

Finally, $\sum_{v\in \{x,y\}} \delta_{A_1}(v)$ is the number of diameters with one end in $\{x, y\}$ and the other end in ${B_1}$, which is, in fact, a vertex in $\partial {B}_1$. Similarly, $\sum_{v\in \{x,y\}} \delta_{B_1}(v)$   is the number of diameters with one end in $\{z, w\}$ and the other in $A^+_1$, which is in fact $\partial A^+_1$. 

We have that the left-hand side of equality \eqref{thecounting} is equals to $2|E(\DDiam_{G_1})|$. Therefore,
$2|E(\DDiam_{G_1})|=2|\tilde{E}(G_1)|$ implying that $|E(\DDiam_{G_1})|=|\tilde{E}(G_1)|$ and, since $|\tilde{E}(G_1)|$ is even (as remarked above) then $|E(\DDiam_{G_1})|$ is also even, as claimed above. Therefore, $G$ cannot have 2-cutting set and therefore $G$ is $3$-connected.

(\emph{Sufficiency}) Suppose that the graph $G$ is 3-connected and simple. Since $V$ is an extremal configuration then, by the (GHS) Theorem \ref{GHS}, $V$ is tight and thus, by Theorem  \ref{2connected}, is a planar graph. Hence, $G$ is a polyhedral graph. Moreover, by Theorem \ref{dualedges}, $G$ admits a canonical involution, and thus $G$ is an involutive polyhedral graph. Therefore, by Lemma \ref{critical}, $\DDiag_G $ is 4-critical. 

 We proceed by contradiction. Let us suppose that $V$ is not strongly critical for the Vázsonyi problem. Then, there is a strongly critical subset  $V_1\subset V$ implying, by the necessity condition, that the 1-skeleton of $\mathcal{B}(V_1)$, say $G_1$, is planar, simple and 3-connected. By the same arguments as above, the latter implies that $G_1$ is an involutive polyhedral graph, and again by Lemma \ref{critical}, $\DDiag_{G_1}$ is 4-critical, contradicting that $\DDiag_G$ is 4-critical.     
\end{proof}

The following result, in terms of Reuleaux polyhedra, implies Conjecture \ref{vazcritical}.

\begin{theorem}\label{Conjballteo}
     Let $V\subset \mathbb{R}^3$ be an extremal set for the Vázsonyi problem. Then, $\mathcal{B}(V)$ is a Reuleaux polyhedron if and only if $V$ is strongly critical.
\end{theorem}

\begin{proof}
Suppose that $V$ is strongly critical. Then, by Lemma \ref{3connected}, the 1-skeleton of $\mathcal{B}(V)$ is simple and $3$-connected and by Theorem \ref{2connected}, is a planar graph. Therefore, by Steinitz's characterization, $\mathcal{B}(V)$ is an standard ball polyhedron. Moreover, since $V$ is an extremal configuration then, by Theorem $\ref{GHS}$, $\verti (\mathcal{B}(V))=V$ implying thus that $\mathcal{B}(V)$ is a Reuleaux polyhedron.

Suppose now that $\mathcal{B}(V)$ is a Reuleaux polyhedron. Then,  $\mathcal{B}(V)$ is a standard ball polyhedron. Since the 1-skeleton of $\mathcal{B}(V)$ has a polytopal structure then, again by Steinitz's characterization, it is simple and $3$-connected, therefore by Lemma \ref{3connected}, $V$ is strongly critical.
\end{proof}

Finally,  we may use Theorem \ref{MeissnerThm} and Theorem \ref{Conjballteo} in order to prove the following corollary.

\begin{corollary}
Let $V\subset \mathbb{R}^3$ be an extremal set for the Vázsonyi problem. Then there is a unique subset of $V$ that is critical for the Vázsonyi problem.
\end{corollary}
\begin{proof}
We proceed by contradiction, suppose there are two different subsets $W_1$ and $W_2$ of $V$ that are critical for the Vázsonyi problem,  then there is a $w_1\in W_1$ such that $w_1\notin W_2$.  Let us observe that $W_2\subset \mathcal{B}(W_1)$ because the diameter of $V$ is 1 and $\mathcal{B}(W_1)$ is a Reuleaux polyhedron by Theorem \ref{Conjballteo}.

Let $r$ be the minimum distance between $w_1$ to all the points in $W_2$ and $\epsilon=\frac{r}{2}$. Then  by Theorem \ref{MeissnerThm},   there is a partition $P_{\epsilon, w_1}=\{P_1,  P_2,  P_3,  P_4\}$ of $\mathcal{B}(W_1)$ (with $\ddiam(P_i)<1$) such that $W_2 \subset P_2\cup P_3\cup P_4$, which is a contradiction because $\DDiam_{W_2}$ is 4-critical by Lemma \ref{critical}.  

\end{proof}

 \subsection{Proof of Theorem \ref{MainTheorem}}\label{MainResults}

We prove our main contribution by analyzing the \emph{minimal} structures for the Borsuk and Vázsonyi problem in $\mathbb{R}^3$, which are astonishingly  the set of vertices of Reuleaux polyhedra in both cases.

\begin{theorem}\label{finiteC}
    Let $V\subset \mathbb{R}^3$ be a finite set of points with $\mid V \mid =n \geq 4$ . The following three statements are equivalent:
    \begin{enumerate}[label=(\roman*)]
        \item\label{1eq} $V$ is strongly critical for the Vázsonyi problem.
        \item\label{2eq} $\DDiam_V$ is 4-critical.
        \item\label{3eq} $\mathcal{B}(V)$ is a Reuleaux polyhedron.
    \end{enumerate}
\end{theorem}
\begin{proof}
    The equivalence $\ref{1eq} \iff \ref{3eq}. $ follows from Theorem \ref{Conjballteo}, and $\ref{3eq} \implies \ref{2eq}$ follows from Lemma \ref{critical}. We shall prove that $\ref{2eq} \implies \ref{1eq}$ \\
    Since $\DDiam_V$ is  4-critical then each $v\in V$ has degree at least 3 in $\DDiam_V$, thus by Theorem \ref{eqtight} $V$ is tight and then we have $V\subset \verti \mathcal{B}(V)$. We consider two cases.
    
        {\bf Case 1)} If $V=\verti \mathcal{B}(V)$, by Theorem \ref{GHS}, $V$ is extremal for the Vázsonyi problem. Suppose that $V$ is not strongly critical for the Vázsonyi problem, then there is a proper subset $V_1$ of $V$, which is strongly critical for the Vázsonyi problem. This implies that $\DDiam_{V_1}$ is 4-critical (since $\ref{1eq}.\implies \ref{2eq}$), contradicting that $\DDiam_V$ is 4-critical.
        
        {\bf Case 2)} Let us prove that  $V\subsetneq \verti \mathcal{B}(V)$ is not possible.  Suppose that $V\subsetneq \verti \mathcal{B}(V)$,  then by Theorem \ref{GHS}, $e(V)<2n-2$. We may assume that $V$ does not have an extremal subset for the Vázsonyi problem, otherwise it would lead a contradiction as in Case 1. 
        
        Let $m_0=(2n-2)-e(V)$ and $v\in \verti \mathcal{B}(V)\backslash V$, then $v$ has to be adjacent to at least 3 diameters (Definition \ref{vertex}), so we can define a new subset $V_1=V\bigcup \{ v\}$ in $\mathbb{R}^3$, having at least 3 more diameters than $V$, so $V_1$ is ``closer'' to become extremal because the difference $m_1:=(2(n+1)-2)-e(V_1)<m_0$. 
        We may repeat this procedure at most $m_0$ times in order to obtain a set $V_t$, with $t\leq m_0$, which has a critical subset $W$ for the Vázsonyi problem (that could be $V_t$ itself).
        
        Let us observe that $W$ is not equal nor a subset of $V$ by the first paragraph of this case,  then there is a $w\in W$ such that $w\notin V$.  Also let us observe that $V\subset V_t \subset \mathcal{B}(V_t)\subset \mathcal{B}(W)$ and $\mathcal{B}(W)$ is a Reuleaux polyhedron because $\ref{1eq} \implies \ref{3eq}$.  
        
        Let $r$ be the minimum distance from $w$ to all the points in $V_t\setminus \{w\}$ and $\epsilon=\frac{r}{2}$, then,  by Theorem \ref{MeissnerThm},  there is a partition $P_{\epsilon, w}=\{P_1,  P_2,  P_3,  P_4\}$ of $\mathcal{B}(W)$ (with $\ddiam(P_i)<1$) such that $V_t\setminus\{w\}\subset P_2\cup P_3\cup P_4$, which is a contradiction because $V\subset V_t\setminus \{w\}$ but $\DDiam_V$ is 4-critical.  

Therefore, $\ref{1eq}$, $ \ref{2eq}$ and $\ref{3eq}$ are equivalent.
\end{proof}

We clearly have that Theorem \ref{MainTheorem} is a straightforward consequence of Theorem \ref{finiteC}.
\subsection{Special configuration of points}\label{SpecialConfig}

Let us consider the following configuration of $8$ points in $\mathbb{R}^3$. Four points, say $w,  x,  y$ and $z$, are the vertices of a regular tetrahedron with edges of length $1$. 

We shall add other appropriate four points, say $a,  b,  c$ and $d$ (this is the tricky part of the construction). 
Let $B(c)$ be the ball with center $c$ and radius 1 and its boundary is the unit sphere $\mathbb{S}(c)$ with center  $c$. The four desired points will lie at the Reuleaux tetrahedron formed by $B(x)\cap B(y) \cap B(w)\cap B(z)$ as follows. Let $p$ (respectively $q$) be the midpoint of the circular-arc edge between $x$ and $y$ (respectively  circular-arc edge between $z$ and $w$).  

\begin{figure}[h]
\centering
\includegraphics[width=.4\linewidth]{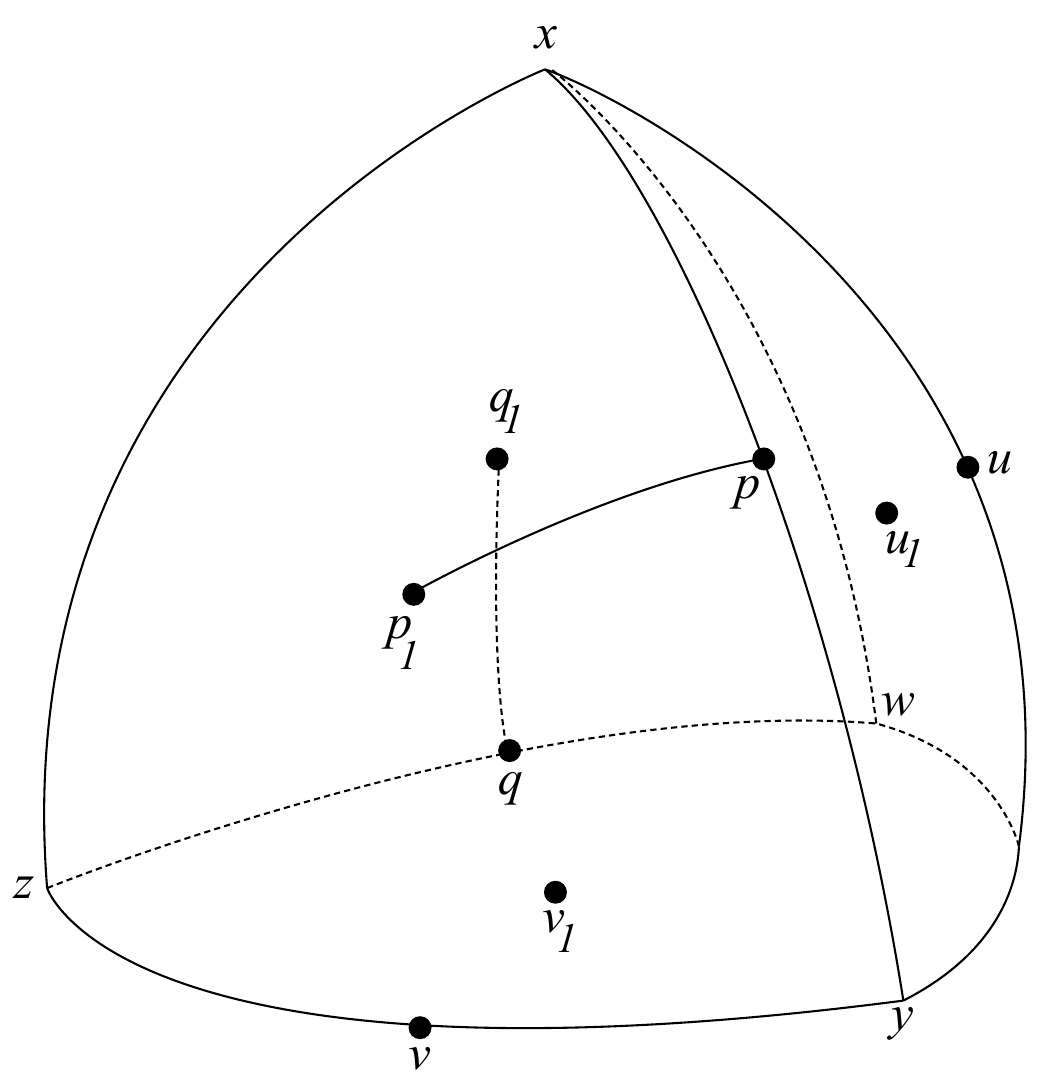}
\caption{Reuleaux Tetrahedron}
\label{FigReul}
\end{figure}

It is known \cite{martini2019bodies} that $\Vert p,q\Vert=\left(\sqrt 3-\frac{\sqrt 2}{2}\right) \approx 1.0249$, see Figure \ref{FigReul}. Let $p_1$ and $q_1$ be points on $\mathbb{S}(w)$ and $\mathbb{S}(y)$ respectively  (in the Reuleaux tetrahedron),  and the Meissner's body of constant width. In such bodies, the segment $[p_1,q_1]$ is not a diameter (see \cite[pp 171-173]{martini2019bodies}), then $\Vert p_1,q_1\Vert <1$.  

Let $A(p,p_1)$ (respectively $A(q,q_1)$) be the circular-arc in $\mathbb{S}(w)$ joining $p$ to $p_1$ (respectively the circular-arc in $\mathbb{S}(y)$ joining $q$ to $q_1$). Let 
$$\begin{array}{lcl}
 \begin{array}{lclc}
\alpha_1: & [0,1] & \longrightarrow & A(p,p_1)\\
 & t & \mapsto &\alpha_1(t)
\end{array}
&
\text{and}
&
\begin{array}{lclc}
\beta_1: & [0,1] & \longrightarrow & A(q,q_1)\\
 & t & \mapsto &\beta_1(t)
\end{array}
\end{array}$$
where $\alpha_1(0) = p, \alpha_1(1)=p_1, \beta_1(0) = q$ and $\beta_1(1)=q_1$.
 
Finally, let 
$$\begin{array}{lclc}
\gamma_1: & [0,1] & \longrightarrow & \mathbb{R} \\
 & t & \mapsto & \gamma_1(t)=\Vert \alpha_1(t),\beta_1(t)\Vert
\end{array}$$

We have that $\gamma_1(t)$ is a continuous function in $[0,1]$. Moreover, since $\gamma_1(0)=\Vert \alpha_1(0),\beta_1(0)\Vert=\Vert p,q\Vert> 1$ and  $\gamma_1(1)=\Vert \alpha_1(1),\beta_1(1)\Vert=\Vert p_1,q_1\Vert< 1$ then, by the Mean Value Theorem, there is $t_1\in [0,1]$ such that $\gamma_1(t_1)=1$. 

We set $a=\alpha_1(t_1)$ and $c=\beta_1(t_1)$.

By using the symmetry of the Reuleaux tetrahedron we construct points $c$ and $d$ in analogous fashion. We consider points $u$ and $v$ (respectively $u_1$ and $v_1$) playing the same role as $p$ and $q$ (respectively $p_1$ and $q_1$), see Figure \ref{FigReul}. 





\begin{figure}[ht]

    \begin{subfigure}[b]{0.4\textwidth}
    \centering
        \definecolor{ffqqqq}{rgb}{1.,0.,0.}
            \begin{tikzpicture}[line cap=round,line join=round,>=triangle 45,x=1.0cm,y=1.0cm, scale=3.5]
            \clip(-0.06307744833005766,-0.06624291357631734) rectangle (1.060094473203906,1.054891767015483);
            \draw [line width=1.pt,color=ffqqqq] (0.,1.)-- (1.,1.);
            \draw [line width=1.pt,color=ffqqqq] (1.,1.)-- (1.,0.);
            \draw [line width=1.pt,color=ffqqqq] (1.,0.)-- (0.,0.);
            \draw [line width=1.pt,color=ffqqqq] (0.,0.)-- (0.,1.);
            \draw [line width=1.pt] (0.5,0.8)-- (0.2,0.3);
            \draw [line width=1.pt] (0.2,0.3)-- (0.8,0.3);
            \draw [line width=1.pt] (0.8,0.3)-- (0.5,0.8);
            \draw [line width=1.pt] (0.5,0.8)-- (0.5,0.5);
            \draw [line width=1.pt] (0.5,0.5)-- (0.2,0.3);
            \draw [line width=1.pt] (0.5,0.5)-- (0.8,0.3);
            \draw [line width=1.pt,color=ffqqqq] (0.5,0.8)-- (1.,1.);
            \draw [line width=1.pt,color=ffqqqq] (0.,1.)-- (0.5,0.5);
            \draw [line width=1.pt,color=ffqqqq] (0.2,0.3)-- (0.,0.);
            \draw [line width=1.pt,color=ffqqqq] (0.8,0.3)-- (1.,0.);
            \draw [color=ffqqqq](0.1096488725292607,0.4129034072270567) node[anchor=north west] {$x$};
            \draw [color=ffqqqq](0.47217979004173605,0.9302292275842829) node[anchor=north west] {$y$};
            \draw [color=ffqqqq](0.5146351005269797,0.5753040004229176) node[anchor=north west] {$z$};
            \draw [color=ffqqqq](0.814967111737407,0.37403438034472625) node[anchor=north west] {$w$};
            \draw (0.9028985297519726,0.9804574361995935) node[anchor=north west] {$a$};
            \draw (0.036619752841273535,0.1141438149234914) node[anchor=north west] {$b$};
            \draw (0.8814501516147551,0.09886107164407401) node[anchor=north west] {$c$};
            \draw (0.050771523003021415,1.018) node[anchor=north west] {$d$};
            \begin{scriptsize}
            \draw [fill=black] (0.,0.) circle (0.8pt);
            \draw [fill=black] (1.,0.) circle (0.8pt);
            \draw [fill=black] (1.,1.) circle (0.8pt);
            \draw [fill=black] (0.,1.) circle (0.8pt);
            \draw [fill=ffqqqq] (0.5,0.8) circle (0.8pt);
            \draw [fill=ffqqqq] (0.2,0.3) circle (0.8pt);
            \draw [fill=ffqqqq] (0.8,0.3) circle (0.8pt);
            \draw [fill=ffqqqq] (0.5,0.5) circle (0.8pt);
            \end{scriptsize}
            \end{tikzpicture}
        \caption{The diameter graph.}
        \label{DiameterGraph}
    \end{subfigure}
    \hfill
        \centering
    \begin{subfigure}[b]{0.3\textwidth}
    \centering
        \includegraphics[width=\textwidth]{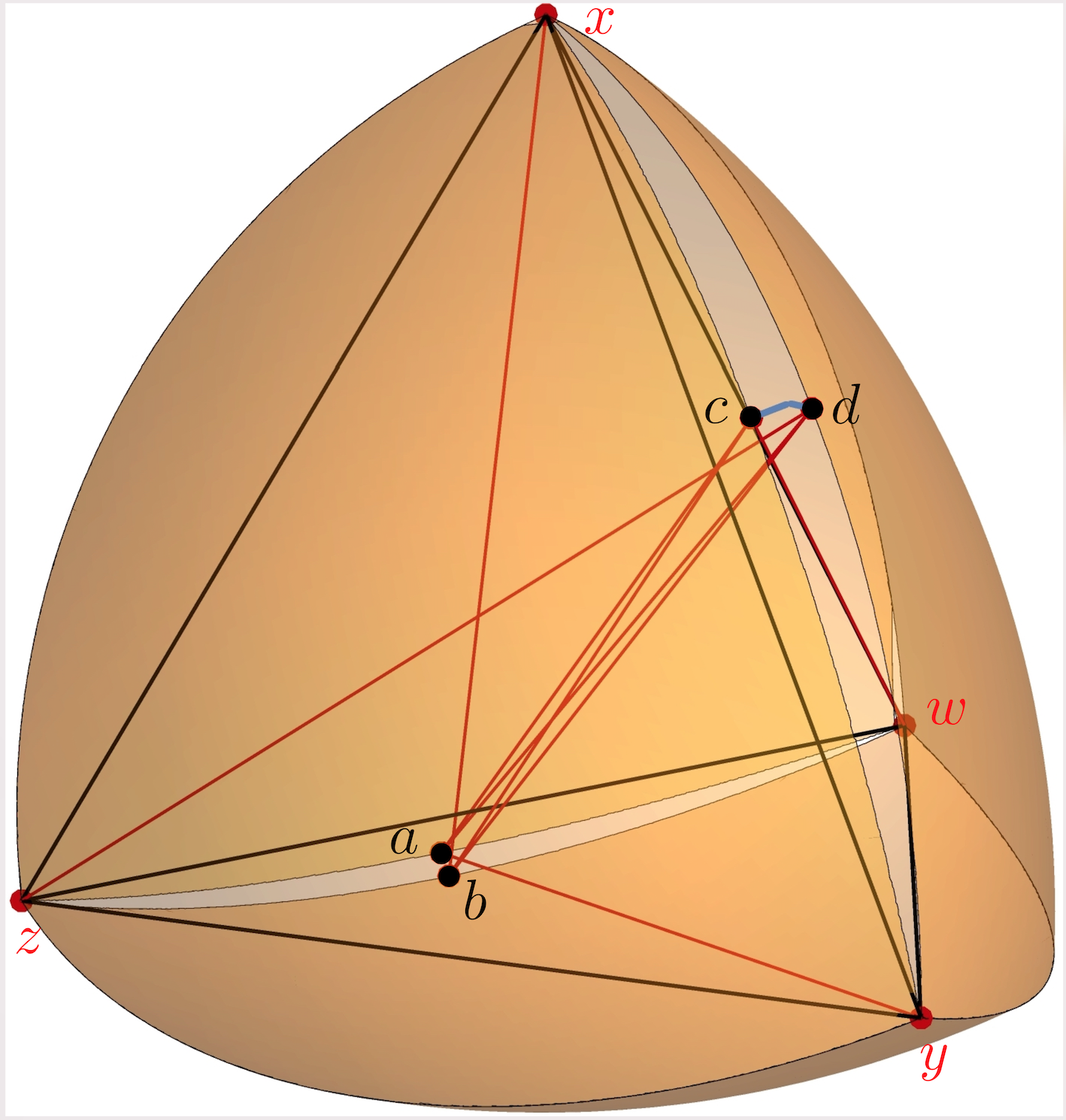}
        \caption{The embedding}
        \label{TheEmbedding}
    \end{subfigure}
     \hfill
    \begin{subfigure}[b]{0.25\textwidth}
     \centering
        \definecolor{qqqqff}{rgb}{0.,0.,1.}
        \definecolor{ffqqqq}{rgb}{1.,0.,0.}
        \begin{tikzpicture}[line cap=round,line join=round,>=triangle 45,x=1.0cm,y=1.0cm, scale=1.1]
        \clip(-0.502740007104047,-0.26086066383508894) rectangle (2.238477940359954,3.8373265464450993);
        \draw [line width=1.pt] (1.,3.5)-- (1.,2.5);
        \draw [line width=1.pt] (0.5,1.)-- (1.5,1.);
        \draw [line width=1.pt] (0.,3.)-- (1.,3.5);
        \draw [line width=1.pt] (1.,3.5)-- (2.,3.);
        \draw [line width=1.pt] (2.,3.)-- (1.,2.5);
        \draw [line width=1.pt] (1.,2.5)-- (0.,3.);
        \draw [line width=1.pt] (1.,2.)-- (0.5,1.);
        \draw [line width=1.pt] (0.5,1.)-- (1.,0.);
        \draw [line width=1.pt] (1.,0.)-- (1.5,1.);
        \draw [line width=1.pt] (1.5,1.)-- (1.,2.);
        \draw [line width=1.pt] (1.,2.)-- (2.,3.);
        \draw [line width=1.pt] (0.,3.)-- (1.,2.);
        \draw [shift={(-1.7056140053522308,2.5685380017840775)},line width=1.pt]  plot[domain=-0.7594138627001588:0.11591275390687458,variable=\t]({1.*3.730648015099667*cos(\t r)+0.*3.730648015099667*sin(\t r)},{0.*3.730648015099667*cos(\t r)+1.*3.730648015099667*sin(\t r)});
        \draw [shift={(3.6987759507617883,2.5662586502539293)},line width=1.pt]  plot[domain=3.024859598502278:3.9018268174705923,variable=\t]({1.*3.724120713995898*cos(\t r)+0.*3.724120713995898*sin(\t r)},{0.*3.724120713995898*cos(\t r)+1.*3.724120713995898*sin(\t r)});
        \draw (1.123476522322596,3.681250989498673) node[anchor=north west] {$a$};
        \draw (0.9801102970404583,2.981422963211218) node[anchor=north west] {$b$};
        \draw (0.1740262807672713,1.3868247078862742) node[anchor=north west] {$c$};
        \draw (1.4579676424743666,1.492432263580917) node[anchor=north west] {$d$};
        \draw [color=ffqqqq](-0.109249727995214018,3.3504052035147587) node[anchor=north west] {$z$};
        \draw [color=ffqqqq](1.8300285398475666,3.3904052035147587) node[anchor=north west] {$w$};
        \draw [color=ffqqqq](0.7645929578654681,1.7759958405558803) node[anchor=north west] {$y$};
        \draw [color=ffqqqq](1.1066538552386682,0.2606420157666748) node[anchor=north west] {$x$};
        \draw [color=qqqqff](0.2578605073011926,3.1997483098969095) node[anchor=north west] {$F_d$};
        \draw [color=qqqqff](1.09851596349636648,3.1904954200344113) node[anchor=north west] {$F_c$};
        \draw [color=qqqqff](0.7255289603690354,2.530766069593369) node[anchor=north west] {$F_x$};
        \draw [color=qqqqff](0.22066005548513357,2.2457167316701625) node[anchor=north west] {$F_w$};
        \draw [color=qqqqff](1.3935751981690093,2.24010917597552) node[anchor=north west] {$F_z$};
        \draw [color=qqqqff](0.8179591831476058,1.527949602755907) node[anchor=north west] {$F_a$};
        \draw [color=qqqqff](0.8111365160636779,0.922333587734504) node[anchor=north west] {$F_b$};
        \draw [color=qqqqff](1.5664734221993633,0.7653220282845107) node[anchor=north west] {$F_y$};
        \begin{scriptsize}
        \draw [fill=ffqqqq] (0.,3.) circle (2.5pt);
        \draw [fill=black] (1.,3.5) circle (2.5pt);
        \draw [fill=black] (1.,2.5) circle (2.5pt);
        \draw [fill=ffqqqq] (2.,3.) circle (2.5pt);
        \draw [fill=black] (0.5,1.) circle (2.5pt);
        \draw [fill=black] (1.5,1.) circle (2.5pt);
        \draw [fill=ffqqqq] (1.,0.) circle (2.5pt);
        \draw [fill=ffqqqq] (1.,2.) circle (2.5pt);
        \end{scriptsize}
        \end{tikzpicture}
        \caption{1-skeleton.}
        \label{skeleton}
    \end{subfigure}
    \caption{Critical configuration of 8 points that is not strongly critical for the Vázsonyi problem.}
    \label{Example8}
   
\end{figure}

Since the original tetrahedron is regular (and each edge is of length one),  the six pairs of points formed by $\{w,x,y,z\}$ are at distance one. Moreover, by construction, $\Vert a,c\Vert=\Vert a,d\Vert=\Vert b,c \Vert=\Vert b,d\Vert=1$. Furthermore,  $\Vert c,w\Vert=\Vert d,z\Vert=\Vert a,y \Vert=\Vert b,x\Vert=1$ since $c\in S(w), d\in S(z), a\in S(y)$ and $b\in S(x)$. It can be checked that the distance of any other pairs of points in $\{a,b,c,d,w,x,y,z\}$ is less than one. The diameter graph is illustrated in Figure \ref{Example8} (a).

The above configuration of 8 points is an extremal Vázsonyi configuration since it contains $(2\times 8) -2=14$ diameters. Moreover, it is critical since all points are adjacent to at least $3$ diameters and there is not dangling edge (see Figure \ref{Example8} center). However, it is not strongly critical since it contains the tetrahedron as an extremal subset.  Moreover, this configuration is an extremal V\'azsonyi configuration but its ball set is not polytopal since it is not 3-connected, for instance $\{z,w\}$ is a 2-cutting set of its 1-skeleton (see Figure \ref{Example8} (c)). The 1-skeleton is indeed planar but just 2-connected.

We computed explicitly the coordinates of the points of such configuration. In order to simplify the calculations, we set the diameter equal to $\sqrt{3}$ and the coordinates for $a,b,c,d$ are approximated with an error of order of $10^{-4}$. 
 $$\begin{array}{ll}
 x=& (0,0, \sqrt 2)\\
 y=& (1,0,0)\\
 w=& (\cos (2\pi/3),\sin(2\pi /3),0)\\
 z=& (\cos (4\pi/3),\sin(4\pi /3),0)\\
 a=& (-0.72849,0,-0.11106)\\
 b=& (-0.68087,0,-0.1784)\\
 c=& (0.7095,-0.03157,0.85524)\\
 d=& (0.7095,0.03157,0.85524)\\
 \end{array}$$


 \section{Concluding remarks}\label{Constructing}

 In this section, we point out some interesting observations and possible  future work on  Reuleaux polyhedra realizations.
 
In  \cite{montejano2020graphs}, the authors computationally proved the validity of Conjecture \ref{ConjReul} up to 14 vertices. They do so by finding first all involutive graphs up to 14 vertices and then constructing explicitly the corresponding desired embedding in each case.  We observe that this list of involutive graphs,  combined with Theorem \ref{finiteC} may allow constructing sets of up to 14 points in $\mathbb{R}^3$ with Borsuk number 4  (extending the examples given in \cite[Lemma 3]{hujter2014multiple} with at most 7 points).

To find the above list of involutive graphs, the authors generated all 3-connected planar graphs and then they searched for the existence of an involutive map in each case. We propose an alternative (more direct) method to find all involutive graphs by using the classification of the family of involutive polyhedra given by Bracho \emph{et al.} \cite[Theorem 6]{bracho2020strongly}.  They showed that if $P$ is an involutive  polyhedral graph,  then there is always an edge $e\in E(P)$ such that $P/\{e\}\setminus\{\tau(e)\}$ is also an involutive polyhedral graph where $\tau$ is the involution and $G\setminus \{f\}$ (respectively denoted by $G/ \{f\}$) denotes the deletion (respectively contraction) of edge $f$ in $G$. The latter implies that any involutive polyhedron can be reduced to a wheel (with an odd number of vertices in the main cycle) by a finite sequence of \emph{delete-contraction} operations (applied simultaneously each time). 

As Tutte \cite{tutte1961theory} remarked, the inverse of the delete-contraction operation corresponds  to diagonalizing faces of the graph and its dual simultaneously. This can be settled as an \emph{add-expansion} operation in $P$ as follows.  

Let $v\in V(P)$ with degree at least 4. Let $F_v$ be the dual face of $v$. Notice that $v$ is a vertex of the dual face $F_w$ for every vertex $w\in F_v$.

 \begin{itemize}
      \item {\bf Split} the vertices of $F_v$ into two paths $P_1$ and $P_2$ with at least 3 vertices each (which is possible since $F_v$ contains at least 4 vertices) with $P_1$ and $P_2$ having only $x$ and $y$ as their common vertices.  {\bf Add} an edge joining $x$ and $y$. Let $F_1$ and $F_2$ be the faces formed by $P_1\cup xy$ and $P_2\cup xy$,  respectively.
    \item {\bf Expand} $v$ into two vertices,  $v_1$ and $v_2$, that is, delete $v$ and add vertices $v_1$ and $v_2$ joined by an edge. Also, for  $i=1, 2$, add an edge joining $v_i$ to a neighbor $w$ of $v$,  such that $\tau(vw)$ (the dual edge of $vw$) is an edge in $P_i$.
 \end{itemize}

We invite the reader to verify that this procedure is the inverse operation of the delete-contraction operation. Let us verify that the resulting graph,  $G'$,  is also an involutive polyhedral graph. We clearly have that $G'$ is a simple, 3-connected, planar graph. Moreover,  the involution $\tau'$ of $G'$ is given by
\[
\tau'(w)=
\begin{cases}
    F_1 &\text{if }w=v_1,\\
      F_2 &\text{if }w=v_2,\\
   \tau(x) \text{ with } v \text{ replaced by the edge } v_1v_2 & \text{if } w=x,\\
   \tau(y) \text{ with } v \text{ replaced by the edge } v_1v_2  & \text{if } w=y,\\
   \tau(w) & \text{otherwise.}
\end{cases}
\]

We thus have that any involutive polyhedral graph can be obtained from an odd wheel by a finite sequence of add-expansion operations. We observe that the latter would lead to a method to construct Reuleaux polyhedra if
Conjecture \ref{ConjReul} were true. Moreover, by Theorem \ref{finiteC}, the former would give infinite families of strongly critical Borsuk configurations as well as strongly critical Vázsonyi configurations.

Also,  by the above, we can deduce that Lemma \ref{critical} gives infinitely many 4-critical graphs that can be actually constructed systematically. It turns out that this infinite family also satisfies the following property that graph theorists might find of interest. We recall that $\delta_{G}(v)$ denotes the degree of vertex $v\in V(G)$,  and $\chi(G)$ is the chromatic number of $G$. \\

\begin{proposition}
    Let $G$ be an involutive polyhedral graph. Then, $\DDiag_G$ is edge 4-critical, that is, it is vertex 4-chromatic and the removal of any edge decreases its chromatic number.
\end{proposition}
\begin{proof} 
We know, by Lemma \ref{critical}, that $\DDiag_G$ is vertex 4-critical. Then, 
$\chi(\DDiag_G)=4$ and $\chi(\DDiag_G\backslash \{v\})<4$ for every $v\in V(\DDiag_G)$. Let $e:=xy\in E(\DDiag_G)$ with $x, y\in V(\DDiag_G)$. We will show that $\chi(\DDiag_G\backslash \{e\})<4$. 

Since $G$ is a polyhedral graph, we have $\delta_{\DDiag_G}(x)\geq 3$ for all $v\in V(\DDiag_G)$. We consider two cases

{\bf Case 1:} $\delta_{\DDiag_G}(x)=3$. Set $F_x:= (y, w_0,w_1)$ and assume the color of $x$ is $c(x)=0$. By Lemma \ref{critical} we know that there is a $3$ coloring of $\DDiag_G\backslash \{x\}$ with colors $\{1,2,3\}$. Suppose $c(y)=1$. If 
$c(w_0), c(w_1)\not= 1$ then we may re-color $x$ with color $c(x)=c(y)=1$ and obtain a proper $3$-coloring of   $\DDiag_G\backslash \{e\}$. If say $c(w_0)=1$ then we may re-color $x$ with color 
$c(x)=j\in \{2,3\}\backslash  c(w_1)$ which yields a proper coloring of $\DDiag_G\backslash \{e\}$. 

{\bf Case 2:} $\delta_{\DDiag_G}(x)\geq 4$. In this case, we apply an add-expansion operation. We do so by expanding $x$ into $v_1$ and $v_2$ in $G$ with $P_1=(w_n,y,w_0)$ and $P_2=(w_0,\dots ,w_n)$ (see the above notation). By the above discussion, the new graph $G'$ is also an involutive polyhedral graph. 

 By construction,  $\DDiag_G\backslash \{e\}$ is a subgraph of $\DDiag_{G'}$. Furthermore, we can obtain $\DDiag_{G'}$ from $\DDiag_G\backslash \{e\}$ by adding a new vertex $z$ and the edges $zw_n$, $zy$ and $zw_0$ (in the above notation, we are taking $v_1=z$ and $v_2=x$).
 
 We thus have that $\DDiag_{G'}\backslash \{z\}=\DDiag_G\backslash \{e\}$. Since $G'$ is also an involutive polyhedral graph,  we know that $\chi(\DDiag_{G'})=4$, and by Lemma \ref{critical}, $\chi(\DDiag_{G'}\backslash\{z\} )=3$, which implies that  $\chi(\DDiag_G\backslash\{e\} )=3$. Therefore,
$\DDiag_{G}$ is edge 4-critical.

\end{proof}
{\bf Acknowledgements}. We would like to thank the referees for many helpful remarks which improved the readability of the paper.

\bibliography{main2}

\end{document}